\documentclass[12pt]{amsart} \textwidth=14.5cm \oddsidemargin=1cm
\usepackage{amsmath}
\usepackage{amsxtra}
\usepackage{amscd}
\usepackage{amsthm}
\usepackage{amsfonts}
\usepackage{amssymb}
\usepackage{eucal}

\input prepictex
\input pictex
\input postpictex

\newtheorem{thm}{Theorem}[section]
\newtheorem{lem}{Lemma}[section]
\newtheorem{cor}{Corollary}[section]
\newtheorem{prop}{Proposition}[section]

\theoremstyle{definition}

\theoremstyle{remark}
\newtheorem{rem}{Remark}[section]

\numberwithin{equation}{section}

\begin{document}

\newcommand{\thmref}[1]{Theorem~\ref{#1}}
\newcommand{\secref}[1]{Section~\ref{#1}}
\newcommand{\lemref}[1]{Lemma~\ref{#1}}
\newcommand{\propref}[1]{Proposition~\ref{#1}}
\newcommand{\corref}[1]{Corollary~\ref{#1}}
\newcommand{\remref}[1]{Remark~\ref{#1}}
\newcommand{\eqnref}[1]{(\ref{#1})}
\newcommand{\exref}[1]{Example~\ref{#1}}

\newcommand{\nc}{\newcommand}
\nc{\on}{\operatorname} \nc{\Z}{{\mathbb Z}} \nc{\C}{{\mathbb C}}
\nc{\oo}{{\mf O}} \nc{\N}{{\mathbb N}} \nc{\bib}{\bibitem}
\nc{\pa}{\partial} \nc{\F}{{\mf F}} \nc{\rarr}{\rightarrow}
\nc{\larr}{\longrightarrow} \nc{\al}{\alpha} \nc{\ri}{\rangle}
\nc{\lef}{\langle} \nc{\W}{{\mc W}} \nc{\gam}{\ol{\gamma}}
\nc{\Q}{\ol{Q}} \nc{\q}{\widetilde{Q}} \nc{\la}{\lambda}
\nc{\ep}{\epsilon} \nc{\g}{\mf g} \nc{\h}{\mf h} \nc{\n}{\mf n}
\nc{\A}{{\mf a}} \nc{\G}{{\mf g}} \nc{\D}{\mc D} \nc{\Li}{{\mc L}}
\nc{\La}{\Lambda} \nc{\is}{{\mathbf i}} \nc{\V}{\mf V}
\nc{\bi}{\bibitem} \nc{\NS}{\mf N}
\nc{\dt}{\mathord{\hbox{${\frac{d}{d t}}$}}} \nc{\E}{\mc E}
\nc{\ba}{\tilde{\pa}}
\def\smapdown#1{\big\downarrow\rlap{$\vcenter{\hbox{$\scriptstyle#1$}}$}}
\nc{\mc}{\mathcal} \nc{\mf}{\mathfrak} \nc{\ol}{\fracline}
\nc{\el}{\ell} \nc{\etabf}{{\bf \eta}} \nc{\x}{{\bf x}}
\nc{\xibf}{{\bf \xi}} \nc{\y}{{\bf y}} \nc{\WW}{\mc W}
\nc{\SW}{\mc S \mc W} \nc{\sd}{\mc S \mc D} \nc{\hsd}{\widehat{\mc
S\mc D}} \nc{\parth}{\partial_{\theta}} \nc{\cwo}{\C[w]^{(1)}}
\nc{\cwe}{\C[w]^{(0)}} \nc{\hf}{\frac{1}{2}}
\nc{\hsdzero}{{}^0\widehat{\sd}} \nc{\hsdpp}{{}^{++}\widehat{\sd}}
\nc{\hsdpm}{{}^{+-}\widehat{\sd}}
\nc{\hsdmp}{{}^{-+}\widehat{\sd}}
\nc{\hsdmm}{{}^{--}\widehat{\sd}} \nc{\gltwo}{{\rm
gl}_{\infty|\infty}} \nc{\btwo}{{B}_{\infty|\infty}}
\nc{\htwo}{{\h}_{\infty|\infty}} \nc{\hglone}{\widehat{\rm
gl}_{\infty}} \nc{\hgltwo}{\widehat{\rm gl}_{\infty|\infty}}
\nc{\hbtwo}{\hat{B}_{\infty|\infty}}
\nc{\hhtwo}{\hat{\h}_{\infty|\infty}} \nc{\glone}{{\rm gl}_\infty}
\nc{\gl}{{\rm gl}} \nc{\ospd}{\mc B} \nc{\hospd}{\widehat{\mc B}}
\nc{\pd}{\mc P} \nc{\hpd}{\widehat{\pd}} \nc{\co}{\mc O}
\nc{\Oe}{\co^{(0)}} \nc{\Oo}{\co^{(1)}} \nc{\sdzero}{{}^0{\sd}}
\nc{\hz}{\hf+\Z} \nc{\vac}{|0 \rangle}

\advance\headheight by 2pt

\title[Lie subalgebras of differential operators on super circle]{Lie
subalgebras of differential operators on the super circle}

\author[Shun-Jen Cheng]{Shun-Jen Cheng$^*$}\thanks{{Mathematics Subject
Classifications (1991)}: 17B67.}
\thanks{$^*$Partially supported by NSC-grant 90-2115-M-006-003 of the R.O.C}
\address{Department of Mathematics, National Taiwan University, Taipei,
Taiwan 106} \email{chengsj@math.ntu.edu.tw}

\author[Weiqiang Wang]{Weiqiang Wang$^{**}$}
\thanks{$^{**}$Partially supported by an NSF grant.}
\address{Department of Mathematics, North Carolina State University,
Raleigh, NC 27695-8205, USA} \email{wqwang@math.ncsu.edu}

\begin{abstract}
We classify anti-involutions of Lie superalgebra $\hsd$ preserving
the principal gradation, where $\hsd$ is the central extension of
the Lie superalgebra of differential operators on the super circle
$S^{1|1}$. We clarify the relations between the corresponding
subalgebras fixed by these anti-involutions and subalgebras of
$\widehat{gl}_{\infty|\infty}$ of types $OSP$ and $P$. We obtain a
criterion for an irreducible highest weight module over these
subalgebras to be quasifinite and construct free field
realizations of a distinguished class of these modules. We further
establish dualities between them and certain finite-dimensional
classical Lie groups on Fock spaces.

\vspace{.3cm}




\end{abstract}

\maketitle \tableofcontents

\section{Introduction}

Extended symmetries have played an important role in conformal
field theories. This leads one to study $\mathcal W$ algebras,
which are higher-spin extensions of the Virasoro algebra, and
their superanalogues, if one takes supersymmetry into account, cf.
\cite{BS, FeF} and the vast literature therein. A fundamental
example of $\WW$ algebras, now known as the $\WW_{1+\infty}$
algebra, appears as the limit, in an appropriate sense, of the
$\WW_N$ algebras, as $N$ goes to $ \infty$ \cite{PRS1}. It has
been realized \cite{PRS2} that the $\WW_{1+\infty}$ algebra can
also be interpreted as the central extension of the Lie algebra of
differential operators on the circle \cite{KP}. Such an
interpretation potentially links $\WW_{1+\infty}$ to geometry in a
direct way, which has yet to be developed. On the other hand, the
super $\WW_{1+\infty}$ algebra, which is a central extension of
the Lie algebra $\sd$ of differential operators on the super
circle $S^{1|1}$ and which will be denoted by $\hsd$ throughout
this paper, is also promising to have connection with mirror
geometry, as it contains the $N=2$ superconformal algebra as a
subalgebra.

The difficulty for the development of a reasonable representation
theory for such Lie (super)algebras lies in the fact that the
graded subspaces are infinite-dimensional in spite of the
existence of a natural principle gradation and thus of a
triangular decomposition. The Cartan subalgebras of such Lie
algebras are also infinite-dimensional. A sensible physical theory
however often requires finite-dimensionality of the graded
subspaces in a representation (which will be referred to as the
{\em quasifiniteness} condition). In \cite{KR1}, Kac and Radul
developed a powerful machinery and initiated a systematic study of
quasifinite representations of $\WW_{1+\infty}$. Subsequently,
there have been further development and extension, cf.
\cite{AFMO1, FKRW, AFMO2, AFMO3, KR2, W2, KWY} and the references
therein. In particular, the Kac-Radul method was further extended
to study the representation theory of the Lie superalgebra $\hsd$
in \cite{AFMO1}.

As explored in \cite{KR1}, the study of representations of
$\WW_{1+\infty}$ is closely related to that of the well-known Lie
algebra $\hglone$. It is well known that $\hglone$ admits natural
subalgebras of type $B,C,D$. In the paper \cite{KWY}, the authors
identified certain distinguished subalgebras of $\WW_{1+\infty}$,
which correspond to these classical subalgebra of $\hglone$, by
first classifying all anti-involutions of $\WW_{1+\infty}$
preserving the principal gradation. Quasifinite representations
and free field realizations of these subalgebras were studied
systematically, and dualities between these subalgebras and
certain finite-dimensional Lie groups were obtained by applying
dualities for classical Lie subalgebras of $\hglone$ in \cite{W1}.
We remark that the free field realization of $\WW_{1+\infty}$ and
duality between $\WW_{1+\infty}$ and finite dimensional general
linear Lie groups were earlier established in \cite{FKRW} and
\cite{KR2}.

The goal of the present paper is to generalize the results of
\cite{KWY} to the $\hsd$ setting. Note that $\hsd$ affords a
canonical principal $\hf \Z$-gradation. More explicitly, we first
classify all anti-involutions of the $\sd$ preserving the
principal gradation. It turns out the theory is much richer in the
super case and we find five families of such anti-involutions.
This is done in Section~\ref{sect_invol}. The anti-involutions
within each family is related to each other by a spectral flow, so
there are essentially five distinct anti-involutions of $\sd$. In
Section~\ref{sect_subalg} we give an explicit description of the
central extension of these five subalgebras fixed by these
anti-involutions, denoted by ${}^0\hsd$ and ${}^{\pm\pm}\hsd$.

By further developing the machinery in \cite{KR1}, we analyze the
structure of parabolic subalgebras of ${}^0\hsd$ and
${}^{\pm\pm}\hsd$, and obtain criterion in terms of certain
differential equations for an irreducible highest weight module of
${}^0\hsd$ and ${}^{\pm\pm}\hsd$ to be quasifinite. This is the
content of Section~\ref{sect_criterion}. Our analysis has a
somewhat different flavor from the one in \cite{KR1}.

The quasifinite representations of ${}^0\hsd$ and
${}^{\pm\pm}\hsd$ turn out to be intimately related to that of the
Lie superalgebra $\hgltwo$ and its subalgebras, $\hospd$ of type
$OSP$ and $\pd_{\pm\pm}$ of type $P$. It is interesting to note
the strange Lie superalgebra of type $P$ (cf. e.g. \cite{K}) makes
a natural appearance here. We develop such a link by exhibiting
certain Lie superalgebra homomorphisms among them in
Section~\ref{sect_embed} (compare \cite{KR1, AFMO1, KWY}). Via
these homomorphisms, the pullback of irreducible quaisifinite
modules of $\hgltwo, \hospd$ and $\pd_{\pm\pm}$ remain irreducible
and quasifinite as modules over ${}^0\hsd$ and ${}^{\pm\pm}\hsd$,
respectively.

In Section~\ref{sect_freefield} we study free field realization of
the Lie superalgebras $\hsd$ and ${}^i\hsd$ in the Fock space
${\mathfrak F}^{\otimes l}$ of $l$ pairs of free fermionic fields
(i.e. $bc$ fields) and $l$ pairs of bosonic ghosts (i.e. $\beta
\gamma$ fields), also cf. \cite{AFMO1}. We further obtain a
duality in the sense of Howe \cite{H1, H2} between $\hsd$ and the
general linear group $GL(l)$. The multiplicity-free decomposition
of ${\mathfrak F}^{\otimes l}$ under the joint action of $\hsd$
and $GL(l)$ is explicitly presented. In particular we obtain
explicit formulas for the joint highest weight vectors in each
isotypic subspace under such a joint action. To achieve this, we
use the technique developed in \cite{CW1} to study in
Section~\ref{sect_howedual} a finite-dimensional duality between a
general linear Lie superalgebra and $GL(l)$ with respect to
certain non-standard Borel subalgebras, which is then adapted to
our infinite-dimensional setting. A similar duality between
${}^0\hsd$ and the Lie group ${\rm Pin}(2l)$ is also obtained. We
remark that the vacuum modules of $\hsd$ and ${}^{+\pm}\hsd$ carry
a natural vertex superalgebra structure.  Finally we conclude this
paper with discussion of some open problems.

{\bf Acknowledgment.} The second author is grateful to the
Max-Planck-Institut f\"ur Mathematik in Bonn and the Center for
Theoretical Sciences in Hsinchu, Taiwan, for providing stimulating
atmosphere and financial support, where part of this work has been
carried out.

{\bf Notation.} $\mathbb N$ stands for the set of positive
intergers, $\Z$ for the set of integers, $\Z_+$ for the set of
non-negative integers, $\Z_2$ for the two-element group ${\Z}/
2\Z$, and $\C$ the field of complex numbers.

\section{Lie superalgebra $\sd$ of differential operators on $S^{1|1}$}
\label{sect_classical}

In this section we set up notation and review the definition of
the Lie superalgebras $\sd$ and $\hsd$, cf. \cite{KR1, AFMO1}.

Let $t$ be an even indeterminate and let
$\partial_t=\frac{d}{dt}$. We denote by ${\D}_{as}$ the
associative algebra of regular differential operators on the
circle $S^1$. It has a linear basis given by
\begin{equation*}
  J^l_k = - t^{l+k} ( \partial_t )^l, \quad
    l \in \Z_{+}, k \in \Z.
\end{equation*}

A different choice of basis of $\D_{as}$ is given by
 $ t^{k} D^l,  l \in \Z_{+}, k \in \Z,$
where $D = t \partial_t$. Note that $f(D) t = t f( D +1)$ for
$f(w) \in \C [w]$ and hence $J^l_k = -t^k [D]_l$, where here and
further we use the notation
\begin{equation}
 [x]_l = x(x-1)\ldots (x-l+1).
   \label{eq_symb}
\end{equation}

Let $\D$ denote the Lie algebra obtained from $\D_{as}$ by taking
the usual bracket of operators $[a, b] =ab - ba.$

Denote by ${\sd}_{as}$ the associative superalgebra of regular
differential operators on the supercircle $S^{1|1}$. Then the
following elements $$
  t^{k+l} (\partial_t)^l \theta \parth,
  t^{k+l} (\partial_t)^l \parth \theta, t^{k+l} (\partial_t)^l \theta,
  t^{k+l} (\partial_t)^l \parth,
  \quad  l \in \Z_+, k \in \Z
$$ form a linear basis of ${\sd}_{as}$, where $\theta$ is an odd
indeterminate. Here of course $\theta$ and $\parth$ commute with
$t$ and $\partial_t$. Note that the odd elements $\theta$ and
$\parth$ generate the four-dimensional Clifford superalgebra with
relation $\theta
\parth+\parth\theta = 1$. Clearly ${\sd}_{as}$ is isomorphic to the tensor
algebra of ${\D}_{as}$ and this Clifford superalgebra.

Denote by $M(1,1)$ the set of all $2 \times 2$ matrices of the
form
\begin{eqnarray}\label{mat_two}
 \left[
  \begin{array}{cc}
   {\alpha}^0   &   {\alpha}^+     \\
   {\alpha}^-   &   {\alpha}^1
 \end{array} \right],
\end{eqnarray}
where ${\alpha}^a  \in \C , a = 0, 1, \pm$, viewed as the
associative superalgebra of linear transformations on the complex
$(1|1)$-dimensional superspace $\C^{1|1}$. Namely, letting $M_a$
$(a = 0, 1, \pm)$ be the matrix of the form \eqnref{mat_two} with
${\alpha}^a = 1$ and $0$ elsewhere, we declare $M_0, M_1$ to be
even and $M_+, M_-$ to be odd elements.  This equips $M(1,1)$ with
a $\Z_2$-gradation. The supertrace $\rm{Str}$ of the matrix
\eqnref{mat_two} is then defined to be ${\alpha}^0 - {\alpha}^1$.

Note that the four-dimensional Clifford superalgebra generated by
$\theta$ and $\parth$ can be identified with ${M}(1,1)$ as
associative superalgebras by making the following identification:
$$
 M_0 = \parth \theta, \;\; M_1 = \theta \parth,\;\; M_+ = \parth, \;\; M_- = \theta.
$$ It follows therefore that one can canonically identify the
associative superalgebra ${\sd}_{as}$ with the associative
superalgebra of $2 \times 2$ (super)matrices with entries in
${\D}_{as}$. By taking the usual super-bracket $[ \;, \; ]$ we
make ${\sd}_{as}$ into a Lie superalgebra, which we denote by
$\sd$. We will adopt the convention that the capital letter $F(D)$
denotes the matrix
\begin{eqnarray}
 \left[
  \begin{array}{cc}
   f_0 (D)  &   f_+  (D)   \\
   f_- (D)  &   f_1  (D)
 \end{array} \right],
   \label{mat_dif}
\end{eqnarray}
where $f_a (w) \in \C [w], a = 0, 1, \pm$, so that we may regard
$F(D)$ as an element of $\sd$.

Denote by $\hsd$ the central extension of $\sd$ by a
one-dimensional vector space with a specified generator $C$. The
commutation relation for homogeneous $F(D)$ and $G(D)$ in $\hsd$
is given by
\begin{eqnarray}
  \left[
     t^r F(D), t^s G(D)
  \right]
    & = & t^{r+s}
    \left(
      F(D + s) G(D) - (-1)^{|F||G|} F(D) G(D +r)
    \right) \nonumber \\
   & + & \Psi
      \left(
        t^r F(D), t^s G(D)
      \right)
      C,
  \label{eq_com}
\end{eqnarray}
where the two-cocycle $ \Psi$ is given by:
\begin{eqnarray}\label{sdcocycle}\nonumber
  \Psi
      \left(
        t^r F(D), t^s G(D)
      \right)
   & = &
   \left\{
      \everymath{\displaystyle}
      \begin{array}{ll}
        \sum_{-r \leq j \leq -1} {\rm Str} \left(
                                            F(j) G(j+r)
                                     \right)
           , & r = -s \geq 0  \\
          0, & r + s \neq 0.
      \end{array}
    \right.  \\
   & = &
       \sum_{-r \leq j \leq -1}
  \left(
        f_0 (j) g_0 (j+r) + f_+ (j) g_- (j+r)
  \right.   \\ \nonumber
   & &
  \left.
       - f_- (j) g_+ (j+r) - f_1 (j) g_1 (j+r)
  \right).
  \label{cocy}
\end{eqnarray}
Here $|F|$ and $|G|$ stand for the degree of $F(D)$ and $G(D)$,
respectively.

We introduce the {\em principal gradation} of $\sd$ and $\hsd$ by
letting the weight of $C$ be 0, the weights of $t^n f(D) \parth
\theta$ and $ t^n f(D) \theta \parth$ be $n$, and the weights of
$t^n f(D) \theta$ and $ t^{n +1} f(D) \parth $ be $n + 1/2$, that
is, the grading is uniquely determined by assigning $wt(C)=0,
wt(D)=0, wt(t)=1, wt(\theta)=1/2$, and $wt (\parth)= -1/2.$ This
equips $\sd$ and $\hsd$ with $\hf\Z$-gradations compatible with
their $\Z_2$-gradations. Namely, if the $\Z_2$-graded
decompositions of $\sd, \hsd$ are given by $$
  \sd = {\sd}_{\bar{0}} \bigoplus {\sd}_{\bar{1}}, \quad
 \hsd = {\hsd}_{\bar{0}} \bigoplus  {\hsd}_{\bar{1}},
$$ then
\begin{eqnarray*}
   {\sd}_{\bar{0}} = \bigoplus_{n \in \Z} {\sd}_n, \quad
   {\sd}_{\bar{1}} = \bigoplus_{n \in \Z + 1/2} {\sd}_n, \\
  {\hsd}_{\bar{0}} = \bigoplus_{n \in \Z} {\hsd}_n, \quad
  {\hsd}_{\bar{1}} = \bigoplus_{n \in \Z + 1/2} {\hsd}_n.
\end{eqnarray*}
Note that each graded subspace $\hsd_n$ is still
infinite-dimensional.

We also have the following triangular decompositions of $\sd$ and
$\hsd$:
\begin{eqnarray*}
   \sd = \sd_{+} \bigoplus \sd_{0}  \bigoplus \sd_{-},  \quad
  \hsd = \hsd_{+} \bigoplus \hsd_{0}  \bigoplus \hsd_{-},
\end{eqnarray*}
where
\begin{eqnarray*}
   \sd_{\pm} = \bigoplus_{j \in \pm \mathbb N /2} \sd_j, \quad                  %
  \hsd_{\pm} = \bigoplus_{j \in \pm \mathbb N /2} {\hsd}_j.
\end{eqnarray*}
\section{Anti-involutions of $\sd$ preserving the principal gradation}
\label{sect_invol}

In this section we classify anti-involutions of ${\sd}_{as}$ which
preserve the principal gradation. An {\em anti-involution}
$\sigma$ of the superalgebra $\sd_{as}$ is an involutive
anti-automorphism of $\sd_{as}$, i.e. $ \sigma^2 = I, \sigma (a X
+ bY) =a\sigma (X) +b \sigma (Y)$ and $ \sigma (XY)= (-1)^{|X||Y|}
\sigma (Y)\sigma (X),$ where $a, b \in \mathbb C, X, Y \in
\sd_{as}$.  We have the following description of anti-involutions
of $\sd_{as}$. The classification here is much more involved than
in the case of $W_{1+\infty}$ as done in \cite{KWY}.

\begin{thm}    \label{th_inv}
Any anti-involution $\sigma$ of $\sd_{as}$ which preserves the
principal gradation is one of the following ($a \in \C$, $b \in
\C^{\times}$):
 \begin{itemize}
  \item[1)] $\sigma_{a,b}( t ) = - t,
       \; \sigma_{a,b}(D) = -D +a+\parth \theta,
       \; \sigma_{a,b}( \theta ) =  b t \parth,
       \; \sigma_{a,b}(\parth) = - (bt)^{-1} \theta$;
  \item[2)] $\sigma_{++,a}( t ) = t,
       \quad \sigma_{++,a}(D) = -D +a,
       \quad \sigma_{++,a}( \theta ) = \theta,
       \quad \sigma_{++,a}(\parth) = - \parth $;
  \item[3)] $\sigma_{+-,a}( t ) = t,
       \quad \sigma_{+-,a}(D) = -D +a,
       \quad \sigma_{+-,a}( \theta ) = - \theta,
       \quad \sigma_{+-,a}(\parth) = \parth $;
  \item[4)] $\sigma_{-+,a}( t ) = -t,
       \quad \sigma_{-+,a}(D) = -D +a,
       \quad \sigma_{-+,a}( \theta ) = \theta,
       \quad \sigma_{-+,a}(\parth) = - \parth $;
  \item[5)] $\sigma_{--,a}( t ) = -t,
       \quad \sigma_{--,a}(D) = -D +a,
       \quad \sigma_{--,a}( \theta ) = - \theta,
       \quad \sigma_{--,a}(\parth) = \parth$.
 \end{itemize}
\end{thm}

The rest of this section is devoted to the proof of
\thmref{th_inv}. The following simple lemma is often used in
computations.

\begin{lem}
Given commuting variables $x, y$ and a polynomial $f$ of one
variable, the following identities hold.
 \begin{eqnarray*}
  f (x + \theta \parth y) & = & f(x) + \theta \parth (f (x +y) - f(x) ), \\
 f ( \parth \theta x + \theta \parth y) & = & \parth \theta f(x) + \theta \parth f(y),  \\
  f (x + \theta \parth y) \parth \theta & = &  f(x) \parth \theta, \\
  f (x + \parth \theta y ) \theta \parth & = &  f(x) \theta \parth.
 \end{eqnarray*}
   \label{lem_ident}
\end{lem}

\begin{proof}
We will prove the first identity. Other identities are obtained
similarly.

If suffices to check for a monomial $f(z) = z^n$.  But in this
case one has
 \begin{align*}
  (x + \theta \parth y)^n
    & = \sum_{k =0}^n {n \choose k} x^k (\theta \parth y)^{n -k}
     = x^n + \sum_{k =1}^n {n \choose k} \theta \parth x^k y^{n -k}   \\
    & = x^n + \theta \parth \left( (x +y)^n - x^n
                       \right),
 \end{align*}
which is the first identity.
\end{proof}

\begin{lem}
Let $\sigma$ be an anti-involution of $\sd_{as}$ preserving the
principal gradation. We have the following two possibilities:
\begin{itemize}
\item[1)]
$\sigma (\theta \parth) = \theta \parth$,
\item[2)] $\sigma (\theta \parth) = \parth \theta $.
   \label{lem_key}
\end{itemize}
\end{lem}

\begin{proof}
 As $\sigma$ preserves the principal gradation, we
 may assume
 \begin{equation}
  \sigma (\theta \parth) = \theta \parth k_0 (D)  + \parth \theta k_1 (D)
  \label{eq_a1}
 \end{equation}
 for some $k_0 (w), k_1 (w) \in \C [w]$. Since $(\theta \parth)^2 = \theta \parth$,
 we have
 \begin{equation}
  \sigma ( \theta \parth )
  = \sigma \left( (\theta \parth)^2 \right)
  = \left( \sigma (\theta \parth)
    \right)^2
  = k_0 (D)^2  \theta \parth  + \parth \theta k_1 (D)^2
  \label{eq_a2}
 \end{equation}
 by using \lemref{lem_ident} and \eqnref{eq_a1}.
 It follows by comparing \eqnref{eq_a1} with \eqnref{eq_a2} that
 \begin{eqnarray}
  k_0 (w) & = & k_0 (w)^2   \label{eq_a3}, \\
  k_1 (w) & = & k_1 (w)^2   \label{eq_a4}.
 \end{eqnarray}
 The solutions of the equations (\ref{eq_a3}) and (\ref{eq_a4}) are
 $ (k_0, k_1) = (0, 1), (1, 0),  (0,0)$ or $ (1, 1)$.
 The first two solutions give rise to what we have listed in 1) and 2).
 The last two solutions
 are easily ruled out since an involution can not send
 $\theta \parth$ to $0$ or $1$.
\end{proof}

Note that the case 1) in \lemref{lem_key} is equivalent to $\sigma
(\parth \theta) = \parth \theta$ while the case 2) in
\lemref{lem_key} is equivalent to $\sigma (\parth \theta) = \theta
\parth $.

In the remainder of this section, we will continue to assume that
$\sigma$ is  an anti-involution of $\sd_{as}$ preserving the
principal gradation. We will analyze the two cases that occurred
in \lemref{lem_key} one by one.
\subsection{The case when $\sigma (\theta \parth) = \theta
\parth$}

This subsection is devoted to the proof of the following.
\begin{prop}
Assume that an anti-involution $\sigma$ of $\sd$ preserving the
principal gradation satisfies $\sigma (\theta \parth) = \theta
\parth$. Then $\sigma$ is one of the $\sigma_{a,b}$ in part 1) of
\thmref{th_inv}.\label{prop_inv1}
\end{prop}

We divide the proof into a series of lemmas.

\begin{lem}
 Assume that $\sigma (\theta \parth) = \theta \parth$. Then
 $ \sigma (\theta ) = b t \parth, \sigma (\parth) = b^{-1} t^{-1}\theta $
 for some nonzero $b \in \C$.
   \label{lem_the}
\end{lem}

\begin{proof}
  We may assume that
 \begin{eqnarray*}
  \sigma ( \theta) & = & \theta h_- (D) + \parth t h_+ (D),\quad \mbox{ for some }
 h_{\pm} (w) \in \C [w],   \\
  \sigma ( \parth) & = & \theta t^{-1} j_- (D) + \parth j_+ (D),
 \quad \mbox{ for some } j_{\pm} (w) \in \C [w].
  \end{eqnarray*}
 Since $[ \theta, \theta \parth ] = - \theta$, we have
 \begin{eqnarray}
 [ \sigma (\theta \parth ), \sigma (\theta ) ] = - \sigma (\theta )
 = - \theta h_- (D) - \parth t h_+ (D).
    \label{eqn_a5}
 \end{eqnarray}
 On the other hand we have
 \begin{eqnarray}
  [ \sigma (\theta \parth ), \sigma (\theta ) ]
  = [ \theta \parth, \theta h_- (D) + \parth t h_+ (D) ]
  = \theta h_- (D) - \parth t h_+ (D).
    \label{eqn_a6}
 \end{eqnarray}
 If follows by comparing \eqnref{eqn_a5} and \eqnref{eqn_a6}
 that $h_- (w) = 0$.
 Similar calculation by using
 $[ \sigma (\theta \parth ), \sigma (\parth ) ] =  \sigma (\parth )$
 shows that $j_+ (w) = 0$. If follows from the following identity
 \begin{eqnarray*}
 \theta \parth
   & = & \sigma (\theta \parth )    \\
   & = & - \sigma (\parth ) \sigma (\theta )   \\
   & = & - \left( \theta t^{-1} j_- (D)
      \right)
      \left( \parth t h_+ (D)
      \right)   \\
   & = & - \theta \parth j_- (D +1) h_+ (D)
 \end{eqnarray*}
 that $j_- (D +1) h_+ (D) = -1$. Since both $ j_- (w)$
 and $h_+ (w)$ are polynomials, $h_+ (w)$ is some nonzero constant $b$
 while $ j_- (w)$ is $ - b^{-1}$.
\end{proof}

\begin{lem}
  Assume that $\sigma (\theta \parth) = \theta \parth$. Then
  $\sigma (t) = - t$.
    \label{lem_t}
\end{lem}

\begin{proof}
  Assume that
 \begin{eqnarray*}
  \sigma (t) & = & \parth \theta t f_0 (D) + \theta \parth t f_1 (D) \quad
  \mbox{for some } f_0 (w), f_1 (w) \in \C [w]   \\
 \sigma (D) & = & \parth \theta  g_0 (D) + \theta \parth  g_1 (D)\quad
  \mbox{for some } g_0 (w), g_1 (w) \in \C [w].
 \end{eqnarray*}
 Then we have
 \begin{eqnarray}
  t & = & \sigma^2 (t)
      = \sigma \left( \parth \theta t f_0 (D) + \theta \parth t f_1 (D)
               \right)    \nonumber \\
    & = & f_0 \left( \parth \theta g_0 (D) + \theta \parth  g_1 (D)
                            \right)
          \left( \parth \theta t f_0 (D) + \theta \parth t f_1 (D)
          \right) \parth \theta \nonumber \\
     & & + f_1 \left( \parth \theta g_0 (D) + \theta \parth   g_1 (D)
               \right)
          \left( \parth \theta t f_0 (D) + \theta \parth t f_1 (D)
          \right) \theta \parth  \nonumber \\
    & = & \parth \theta f_0 \left(  g_0 (D)
                         \right)
                         \left( t f_0 (D)
                         \right)
        + \theta \parth  f_1 \left(  g_1 (D)
                         \right)
                         \left( t f_1 (D)
                         \right)   \nonumber \\
    & = & \parth \theta t f_0 \left(  g_0 (D +1)
                         \right)
                         \left( f_0 (D)
                         \right)
        + \theta \parth t f_1 \left(  g_1 (D +1)
                         \right)
                         \left( f_1 (D)
                         \right) .
   \label{eq_t1}
 \end{eqnarray}
 As we may rewrite $ t = \parth \theta  t + \theta \parth t$,
 it follows from \eqnref{eq_t1} that
 $$
   f_0 \left(  g_0 (w +1)  \right) \left( f_0 (w) \right)
    =  f_1 \left(  g_1 (w +1)  \right) \left( f_1 (w) \right) = 1.
 $$
 This implies that $(f_0 (w), f_1 (w) ) = (1, 1), (-1, -1), (1, -1)$
 or $(-1, 1)$. In the first two cases
 we get $\sigma (t) = \pm t$.
 The last two solutions for $(f_0 (w), f_1 (w) )$ give rise to
 $\sigma (t) = \pm ( \parth \theta t - \theta \parth t )$. We first
 claim that these last
 two cases cannot happen. For example let us take
 $\sigma (t) =  \theta \parth t - \parth \theta t $ (the case
 $\sigma (t) = - ( \theta \parth t - \parth \theta t )$ is similar).
 On one hand by \lemref{lem_the} we have
 \begin{eqnarray*}
   \sigma (t) \sigma ( \theta)
   &= & \left( \theta \parth t - \parth \theta t \right)
        \left( b t \parth \right) = -b t^2 \parth,   \\
   \sigma ( \theta) \sigma (t)
   &= & \left( b t \parth \right)
        \left( \theta \parth t - \parth \theta t \right)
  = b t^2 \parth.
 \end{eqnarray*}
 So $ \sigma (t) \sigma ( \theta) = - \sigma ( \theta) \sigma (t)$.
 On the other hand,
 $ \sigma (t) \sigma ( \theta)
   = \sigma ( \theta) \sigma (t)$ since $ \theta t = t \theta$. This is a contradiction.

 We claim that $\sigma (t) =  t$ cannot occur either.
 Indeed if $\sigma (t) =  t$ then by \lemref{lem_the}
 \begin{eqnarray*}
  \theta = \sigma^2 (\theta) = \sigma ( b t \parth)
      = b \sigma ( \parth)\sigma (t)
      = b ( - b^{-1} t^{-1} \theta) t
      = - \theta
 \end{eqnarray*}
  which is a contradiction.
\end{proof}

\begin{lem}
  Assume that $\sigma (\theta \parth) = \theta \parth$. Then
  $\sigma (D) = -D + a + \parth \theta$ for some $a \in \C$.
    \label{lem_D}
\end{lem}

\begin{proof}
  Assume that $\sigma (D) = \parth \theta g_0 (D) + \theta \parth  g_1 (D)$
 for some $g_0 (w), g_1 (w) \in \C [w]$. Then we have
 \begin{eqnarray*}
   D
   & = & \sigma^2 (D) \\
   & = & \sigma (\parth \theta g_0 (D) + \theta \parth  g_1 (D))  \\
   & = & g_0 (\parth \theta g_0 (D) +  \theta \parth  g_1 (D) ) \parth \theta
   + g_1 ( \parth \theta g_0 (D) + \theta \parth  g_1 (D) ) \theta \parth \\
   & = & \parth \theta g_0 ( g_0 (D) ) + \theta \parth g_1 ( g_1 (D) )  \\
   & = & g_0 ( g_0 (D) ) +  \theta \parth ( g_1 ( g_1 (D) ) - g_0 ( g_0 (D) ) ).
 \end{eqnarray*}
 Thus $g_0 ( g_0 (w) ) = g_1 ( g_1 (w) ) = w$.
 This implies $g_1 (w) = - w +a$ for some $a \in \C$
 or $g_1 (w) = w$.

 We have $\sigma (\theta) \sigma (D) = \sigma (D) \sigma (\theta)$
 as $D \theta = \theta D$.  On the other hand, we can compute directly that

 \begin{eqnarray*}
 \sigma (\theta) \sigma (D)
   & = &  (b t \parth) \left( \parth \theta  g_0 (D) + \theta \parth g_1 (D)
                  \right)
     =  b \parth t g_1 (D), \\
 \sigma (D) \sigma (\theta)
   & = &  \left( \parth \theta g_0 (D) + \theta \parth  g_1 (D)
          \right)
          (b t \parth )
     =  b \parth t g_0 (D +1).
 \end{eqnarray*}
 It follows that $ g_0 (w +1) = g_1 (w)$.
 In the case that $g_1 (w) = - w + a$, $g_0 (w) = g_1 (w -1) = -w +a +1$.
 Thus $ \sigma (D) = -D +a +\parth \theta$. We claim that
 $g_1 (w) = w$ cannot occur. Indeed if $g_1 (w) = w$, then
 $g_0 (w) = w -1$ and $\sigma (D) = D - \parth \theta$.
 Together with \lemref{lem_t} and the fact that
 $[t, D] = -t$ this implies that
 $$ [D, - t] = [ \sigma (D), \sigma (t) ] = - \sigma (t) = t.$$
This contradicts the fact that $[t, D] = -t$.
\end{proof}

Proposition~\ref{prop_inv1} now follows from \lemref{lem_the},
\lemref{lem_t} and \lemref{lem_D}.
\subsection{The case when $\sigma (\theta \parth) = \parth
\theta$}

This subsection is devoted to the proof of the following.

\begin{prop}
Assume that an anti-involution $\sigma$ of $\sd$ preserving the
principal gradation satisfies  $\sigma (\theta \parth) = \parth
\theta$. Then
  $\sigma$ is either $\sigma_{++,a}$, $\sigma_{+-,a}$,
 $\sigma_{-+,a}$ or $\sigma_{--,a}$ in \thmref{th_inv}.
   \label{prop_inv2}
\end{prop}

We divide the proof into a series of lemmas.
\begin{lem}
 Assume that $\sigma (\theta \parth) = \parth \theta$. Then
 $ \sigma (\theta ) = \pm \theta$ and $ \sigma (\parth) = \mp \parth $.
   \label{lem_the2}
\end{lem}

\begin{proof}
   Assume that
 \begin{eqnarray*}
 \sigma ( \theta) & = & \theta h_- (D) + \parth t h_+ (D), \mbox{ for some }
 h_+ (w), h_- (w) \in \C [w],\\
 \sigma ( \parth) & = & \theta t^{-1} j_- (D) + \parth j_+ (D),
 \mbox{ for some } j_+ (w), j_- (w) \in \C [w].
 \end{eqnarray*}
 Since $[ \theta, \theta \parth ] = - \theta$, we have
 \begin{eqnarray}
 [ \sigma (\theta \parth ), \sigma (\theta ) ] = - \sigma (\theta )
 = - \theta h_- (D) - \parth t h_+ (D).
    \label{eqn_b5}
 \end{eqnarray}
 On the other hand, we also have
 \begin{eqnarray}
  [ \sigma (\theta \parth ), \sigma (\theta ) ]
  = [ \parth \theta, \theta h_- (D) + \parth t h_+ (D) ]
  = - \theta h_- (D) + \parth t h_+ (D).
    \label{eqn_b6}
 \end{eqnarray}
 If follows by comparing \eqnref{eqn_b5} and \eqnref{eqn_b6}
 that $h_+ (w) = 0$.
 Similar calculation by using
 $[ \sigma (\theta \parth ), \sigma (\parth ) ] =  \sigma (\parth )$
 shows that $j_- (w) = 0$. If follows from
 \begin{eqnarray*}
  \parth \theta
  = \sigma (\theta \parth )
  = - \sigma (\parth ) \sigma (\theta )
  = - \left( \parth  j_+ (D)
      \right)
      \left( \theta h_- (D)
      \right)
  = - \parth \theta j_+ (D ) h_- (D)
 \end{eqnarray*}
 that $j_+ (w ) h_- (w) = -1$. Since both $ j_+ (w)$
 and $h_- (w)$ are polynomials, $h_- (w)$ is some nonzero constant $b$
 while $ j_+ (w)$ is $ - b^{-1}$. Noting that
 $\theta = \sigma^2 (\theta) = b^2 \theta$, we have $ b = \pm 1$.
\end{proof}

\begin{lem}
  Assume that $\sigma (\theta \parth) = \parth \theta$. Then
  $\sigma (t) = \pm t$.
    \label{lem_t2}
\end{lem}

\begin{proof}
  Assume that
 \begin{eqnarray*}
 \sigma (t) & = & \parth \theta t f_0 (D) + \theta \parth t f_1 (D)
 \mbox{ for some } f_0 (w), f_1 (w) \in \C [w],\\
 \sigma (D) & = & \parth \theta t g_0 (D) + \theta \parth  g_1 (D)
 \mbox{ for some } g_0 (w), g_1 (w) \in \C [w].
 \end{eqnarray*}
 We calculate
 \begin{eqnarray*}
  t & = & \sigma^2 (t)  \\
    & = & \sigma \left( \parth \theta t f_0 (D) + \theta \parth t f_1 (D)
                 \right) \\
    & = & f_0 \left( \parth \theta g_0 (D) + \theta \parth g_1 (D)
                            \right)
          \left( \parth \theta  t f_0 (D) + \theta \parth t f_1 (D)
          \right) \theta \parth  \\
     & & + f_1 \left( \parth \theta   g_0 (D) + \theta \parth   g_1 (D)
               \right)
          \left( \parth \theta t f_0 (D) + \theta \parth t f_1 (D)
          \right) \parth \theta \\
    & = & \theta \parth  f_0 (  g_1 (D) )
                         ( t f_1 (D) )
        + \parth \theta f_1 (  g_0 (D) )
                         ( t f_0 (D) )  \\
    & = & \theta \parth t f_0 (  g_1 (D +1) )
                          ( f_1 (D) )
        + \parth \theta t f_1 (  g_0 (D +1) )
                           ( f_0 (D) ).
 \end{eqnarray*}
 As we may rewrite $t = \parth \theta t +  \theta \parth t $, it follows that
 $$
   f_0 \left(  g_1 (w +1)  \right)  f_1 (w)
    =  f_1 \left(  g_0 (w +1)  \right) f_0 (w)  = 1.
 $$
 This implies that $f_0 (w)$ and $ f_1 (w)$ are constants. Let
 us put $f_0  = c$ $ (c \in \C)$ then $f_1 = c^{-1}$ and so
 $\sigma (t) = c \theta \parth t + c^{-1} \parth \theta t$.

 It follows from $t \theta = \theta t$ that
 $\sigma (t) \sigma (\theta) = \sigma (\theta) \sigma (t)$. On the other hand
 we have by \lemref{lem_the2}
 \begin{eqnarray*}
  \sigma (t) \sigma (\theta)
   & = & (c \theta \parth + c^{-1} \parth \theta) (\pm \theta)
     =   \pm c \theta,\\
  \sigma (\theta) \sigma (t)
   & = & (\pm \theta) (c \theta \parth + c^{-1} \parth \theta)
     =   \pm c^{-1} \theta .
 \end{eqnarray*}

 This implies $c = c^{-1}$ i.e. $c = \pm 1$. Thus
 $\sigma (t) = c \theta \parth t + c^{-1} \parth \theta t = \pm t.$
\end{proof}

\begin{lem}
  Assume that $\sigma (\theta \parth) = \parth \theta $. Then
  $\sigma (D) = -D + a $ for some $a \in \C$.
    \label{lem_D2}
\end{lem}

\begin{proof}
  Assume that $\sigma (D) = \parth \theta g_0 (D) + \theta \parth  g_1 (D)$
 for some $g_0 (w), g_1 (w) \in \C [w]$. We have
 \begin{eqnarray*}
   D
   & = & \sigma^2 (D) \\
   & = & \sigma (\parth \theta  g_0 (D) + \theta \parth  g_1 (D))  \\
   & = & g_0 ( \parth \theta  g_0 (D) + \theta \parth g_1 (D) )  \theta \parth
   + g_1 ( \parth \theta  g_0 (D) + \theta \parth  g_1 (D) ) \parth \theta \\
   & = & \theta \parth g_0 ( g_1 (D) ) + \parth \theta g_1 ( g_0 (D) )  \\
   & = & g_0 ( g_1 (D) ) + \parth \theta ( g_1 ( g_0 (D) ) - g_0 ( g_1 (D) ) ).
 \end{eqnarray*}
 It follows that
 \begin{eqnarray*}
  g_0 ( g_1 (w) ) = g_1 ( g_0 (w) ) = w.
  \label{eqn_ex1}
 \end{eqnarray*}
 This implies that $g_0 ( g_0 (w)) = w$ and thus
 $g_0 (w) = - w +a$ for some $a \in \C$ or $g_0 (w) = w$.

 From $D \theta = \theta D$ we conclude that
 $\sigma (\theta) \sigma (D) = \sigma (D) \sigma (\theta)$.
 We calculate
 \begin{eqnarray*}
 \sigma (\theta) \sigma (D)
   & = &  (\pm \theta) \left( \theta \parth  g_0 (D) + \parth \theta g_1 (D)
                  \right)
     =  \pm \theta  g_1 (D),\\
 \sigma (D) \sigma (\theta)
   & = &  \left( \theta \parth  g_0 (D) + \parth \theta g_1 (D)
          \right)
          (\pm \theta )
     =  \pm \theta g_0 (D).
 \end{eqnarray*}
 Thus $ g_0 (w) = g_1 (w)$. We claim that $g_0 (w) = w$ is impossible.
 Indeed if $g_0 (w) = w$, then $\sigma (D) = D$. So by \lemref{lem_t2}
 (to be definite we choose $\sigma (t) = t$, the case $\sigma (t) = -t$
 is the same) it follows that
 $$
  -t =  - \sigma (t) =  [ \sigma (D), \sigma (t) ] = [D, t] = t.
 $$
 which is a contradiction.
 Thus $g_0 (w) = g_1 (w) = - w +a$ and so $ \sigma (w) = -w +a$.
\end{proof}

Now Proposition~\ref{prop_inv2} follows from \lemref{lem_the2},
\lemref{lem_t2} and \lemref{lem_D2}.

It is straightforward to check that the $\sigma$'s given in
\thmref{th_inv} are indeed anti-involutions of $\sd_{as}$.
Combining \lemref{lem_key} with \propref{prop_inv1} and
\propref{prop_inv2}, we have proved \thmref{th_inv}.
\section{Subalgebras of $\sd$ fixed by anti-involutions}
\label{sect_subalg}

Given an anti-involution $\sigma$ of $\sd$, one can check easily
that $$
  \sd^{\sigma} \equiv \{ e \in \sd| \sigma (e) = -e \}
$$ is a Lie subalgebra of $\sd$. If $\sigma$ preserves the
principal $\hf\Z$-gradation of $\sd$, then $\sd^{\sigma} $
inherits the principal $\hf \Z$-gradation from $\sd$: $$
  \sd^{\sigma} = \oplus_{n \in \hf \Z} \sd^{\sigma}_n.
$$ Denote by $\cwo$ the space of all odd polynomials in $ \C [w]$,
and by $\cwe$ the space of all even polynomials in $ \C [w]$. Let
$ \overline{k} = 0$ if $k$ is an odd integer and $ \overline{k} =
1$ if $k$ is even.  The purpose of this section is to give an
explicit description for the Lie superalgebras $\sd^{\sigma}$,
where $\sigma$ is an involution given in \thmref{th_inv}. We also
provide a set of linear basis elements for these superalgebras,
which are canonical in connections to vertex algebras and free
field realizations.

\begin{prop}\label{sdzeroform}
 We have the following description for the graded subspaces
 $\sd_j^{\sigma_{a,b} }$ and $\sd_{j - 1/2}^{\sigma_{a,b} }$
 $ (j \in \Z)$ of \,$\sd^{\sigma_{a,b} }$:
 \begin{align*}
  &\sd^{\sigma_{a,b} }_j
   = \{ \parth \theta t^j f( D + (j -a -1)/2 ) +
       \theta \parth  t^j g( D + (j -a)/2 )
  \mid f , g \in \C^{(\overline{j})} [w]  \} ,         \\
  &\sd^{\sigma_{a,b} }_{j - 1/2}
   = \{ \theta t^{j -1} g_0 (D) + \parth (-t)^j b g_0 ( -D -j +a +1),
       g_0 \in \C[w]   \}.
 \end{align*}
\end{prop}

\begin{proof}
 Assume that
 $ A = \parth \theta t^j f_0 (D) + \theta \parth t^j f_1 (D)$ and
 $ B = \theta t^{j -1} g_0 (D) + \parth t^j g_1 (D)$
 belong to $ \sd^{\sigma_{a,b}}$.
 \begin{eqnarray*}
  \sigma_{a,b} (A)
 & = & f_0 ( -D +a + \parth \theta) (-t)^j \parth \theta
     + f_1 ( -D +a + \parth \theta) (-t)^j  \theta \parth  \\
 & = & (-t)^j f_0 ( -D -j +a + \parth \theta) \parth \theta
     + (-t)^j f_1 ( -D -j +a + \parth \theta) \theta \parth  \\
 & = & (-t)^j f_0 ( -D -j +a + 1) \parth \theta
     + (-t)^j f_1 ( -D -j +a) \theta \parth.
\end{eqnarray*}

 It follows by comparing with $\sigma_{a,b} (A) = - A$ that
 \begin{eqnarray*}
   f_0 (w) & = & (-1)^{j -1} f_0 ( -w -j +a + 1), \\
  f_1 (w)  & = & (-1)^{j -1} f_1 ( -w -j +a).
 \end{eqnarray*}

 Define $f (w) = f_0 (w - (j -a -1)/2 )$ and
 $g (w) = f_1 (w - (j -a)/2 )$. We see that
 $f (w), g (w) \in  \C^{(\overline{j})} [w] $.
 On the other hand, we calculate that
 \begin{eqnarray*}
  \sigma_{a,b} (B)
  & = & g_0 ( -D +a + \parth \theta) ( -t)^{j -1} (bt \parth)
      + g_1 ( -D +a + \parth \theta) ( -t)^j ( - b^{-1} t^{-1} \theta) \\
  & = & (-1)^{j -1} b g_0 ( -D +a +1) t^j \parth
      + (-1)^{j -1} b^{-1} g_1 ( -D +a) t^{j -1} \theta \\
  & = &  (-1)^{j -1} b t^j g_0 ( -D -j +a +1) \parth \\
     && + (-1)^{j -1} b^{-1} t^{j -1} g_1 ( -D -j +a +1) \theta.
 \end{eqnarray*}

 By comparing with $\sigma_{a,b} (B) = - B$ we have
 \begin{eqnarray*}
  g_0 (D) & = & (-1)^{j} b^{-1} g_1 ( -D -j +a +1), \\
  g_1 (D) & = & (-1)^{j} b g_0 ( -D -j +a +1).
 \end{eqnarray*}
 We observe that these two equations are equivalent to each other.
\end{proof}

Letting $ t \mapsto t, \partial_t \mapsto \partial_t, \theta
\mapsto \alpha \theta$ and $ \parth \mapsto \alpha \parth$ defines
an automorphism $\sharp_{\alpha} $ $(\alpha \in \C^{\times})$ of
$\sd$. It is easy to verify that $$
  \sigma_{a,b} \sharp_{\alpha}
    = \sharp_{\alpha^{-1} } \sigma_{a,b} = \sigma_{a,b \alpha}.
$$ Denote by $\Theta_{s}$ the automorphism of $\sd$ defined by
letting $ t \mapsto t, D \mapsto D + s, \theta \mapsto \theta$ and
$
\parth \mapsto \parth$. Clearly
\begin{eqnarray*}
  \sigma_{a, b} \cdot \Theta_{s} = \sigma_{a+ s, b}
   = \Theta_{-s} \cdot \sigma_{a, b}.
\end{eqnarray*}
In this way, we may regard the $\sigma_{a, b}$ for different $a,b$
are related to each other by some spectral flow (cf. e.g.
\cite{AFMO1}). It follows that the Lie superalgebra
$\sd^{\sigma_{a,b}}$ is isomorphic to $\sd^{\sigma_{a',b'}}$ for
all $a, a' \in \C, b, b' \in \C^{\times}$. It turns out that a
convenient choice is putting $a= -1, b= 1$, which we will fix from
now on and will denote $\sd^{\sigma_{-1, 1}}$ by $\sdzero$. Also
the $2$-cocycle \eqnref{sdcocycle} of $\sd$ restricted to
$\sdzero$ gives rise to a central extension of $\sdzero$, which we
will denote by $\hsdzero$.

We introduce a canonical spanning set of $\sdzero$ that will be
used later. Define
\begin{eqnarray*}
 {}^0 W^{n}_{0,k} & = & t^k \left(
                    [D]_n + (- 1)^{k +1} [ -D -k]_n
                \right) \parth \theta,  \\
 {}^0 W^{n}_{1,k} & = & t^k \left(
                    [D]_n + (- 1)^{k +1} [ -D -k -1]_n
                \right) \theta\parth,  \\
 {}^0 W^{n}_{{\times},k} & = & \theta t^{k} [D]_n
                     + (-1)^{k+1} \parth t^{k+1}  [ -D -k - 1]_n.
   \label{eq_basisw}
\end{eqnarray*}
Then ${}^0 W^{n}_{a,k}$ $(k \in \Z, n \in \mathbb N, a = 0, 1,
{\times})$ span $\sdzero$.

\begin{prop}\label{zeroform1}
 We have the following description for the graded subspaces
$\sd_j^{\sigma_{++,a} }$ and $\sd_{j-1/2}^{\sigma_{++,a} }$ $ (j
\in \Z)$ of \, $\sd^{\sigma_{++,a} }$:
\begin{eqnarray*}
 \sd^{\sigma_{++,a} }_j
  & = & \left \{
                 \parth \theta t^j f ( D) - \theta \parth t^j f(-D -j +a)
                 \mid f (w) \in \C [w]
        \right \} ,      \\
  \sd^{\sigma_{++,a} }_{j - 1/2}
  & = & \left \{ \theta t^{j -1} f (D + (j -a -1)/2)
                     + \parth t^j g (D + (j -a)/2) \right. \\
   && \quad\quad\quad\quad\quad\quad\quad\quad\quad\quad\quad\quad\quad
  \left.  \mid  f(w)\in \C^{(1)}[w], g (w)  \in \C^{(0)} [w]
         \right \}.
 \end{eqnarray*}
\end{prop}

\begin{proof}
 Assume that
 $ A = \parth \theta t^j f (D) + \theta \parth t^j f_1 (D)$ and
 $ B = \theta t^{j -1} g_- (D) + \parth t^j g_+ (D)$ belong
 to $ \sd^{\sigma_{++,a}}$.
 \begin{eqnarray*}
  \sigma_{++,a} (A)
 & = & \theta \parth f ( -D +a) t^j
     + \parth \theta f_1 ( -D +a) t^j    \\
 & = & \theta \parth t^j f ( -D -j +a)
      + \parth \theta t^j f_1 ( -D -j +a).
 \end{eqnarray*}
 We obtain by comparing with $ \sigma_{++,a} (A) = -A$ that
 \begin{eqnarray*}
  - f (D)  =  f_1 ( -D -j +a), \quad
  - f_1 (D) =  f ( -D -j +a).
 \end{eqnarray*}
 The two equations are clearly equivalent.

 On the other hand, we have
 \begin{eqnarray*}
  \sigma_{++,a} (B)
 & = & g_- ( -D +a) t^{j -1} \theta - g_+ ( -D +a) t^j \parth   \\
 & = & \theta t^{j -1} g_- ( -D -j +1 +a) + \parth t^j g_+ ( -D -j +a).
 \end{eqnarray*}
 We obtain by comparing with $ \sigma_{++,a} (B) = -B$ that
 \begin{eqnarray*}
  - g_- (D) & = & g_- ( -D -j +1 +a), \\
  g_+ (D) & = & g_+ ( -D -j +a).
 \end{eqnarray*}
 If we let $f (w) = g_- (w + (a +1 -j)/2 )$ and $g(w) = g_+ (w + (a -j)/2 )$,
 we see that $ f (w)\in \C^{(1)} [w], g(w) \in \C^{(0)} [w]$.
\end{proof}

In a similar way we determine $\sd^{\sigma_{+-,a} }$ as follows.
We omit the arguments, as they are parallel to the ones given
above.
\begin{prop}\label{zeroform2}
 We have the following description for the graded subspaces
 $\sd_j^{\sigma_{+-,a} }$ and $\sd^{\sigma_{+-,a} }_{j - 1/2}$
 $ (j \in \Z)$ of \,$\sd^{\sigma_{+-,a} }$:
  \begin{eqnarray*}
 \sd^{\sigma_{+-,a} }_j
  & = & \left \{
                 \parth \theta t^j f ( D) - \theta \parth t^j f(-D -j +a)
                 \mid f (w) \in \C [w]
        \right \} ,      \\
  \sd^{\sigma_{+-,a} }_{j - 1/2}
  & = & \left \{ \theta t^{j -1} f (D + (j -a -1)/2)
                     + \parth t^j g (D + (j -a)/2) \right. \\
   && \quad\quad\quad\quad\quad\quad\quad\quad\quad\quad\quad\quad\quad
  \left.  \mid  f(w ) \in \C^{(0)}[w] , g (w)  \in \C^{(1)} [w]
         \right \}.
 \end{eqnarray*}
\end{prop}

\begin{prop}\label{zeroform3}
 We have the following description for the graded subspaces
 $\sd_j^{\sigma_{-+,a} }$ and $ \sd^{\sigma_{-+,a} }_{j - 1/2}$
 $(j \in \Z)$ of \,$\sd^{\sigma_{-+,a} }$:
  \begin{eqnarray*}
 \sd^{\sigma_{-+,a} }_j
  & = & \left \{
                 \parth \theta t^j f ( D) - ( -1)^j \theta \parth t^j f(-D -j +a)
                 \mid f (w) \in \C [w]
        \right \} ,      \\
  \sd^{\sigma_{-+,a} }_{j - 1/2}
  & = & \left \{ \theta t^{j -1} f (D + (j -a -1)/2)
                     + \parth t^j g (D + (j -a)/2) \right. \\
   && \quad\quad\quad\quad\quad\quad\quad\quad\quad\quad\quad\quad\quad
  \left.  \mid  f(w), g (w)  \in \C^{(\overline{j +1})} [w]
         \right \}.
 \end{eqnarray*}
\end{prop}
\begin{proof}
 Assume that
 $ A = \parth \theta t^j f (D) + \theta \parth t^j f_1 (D)$ and
 $ B = \theta t^{j -1} g_- (D) + \parth t^j g_+ (D)$ belong
 to $ \sd^{\sigma_{-+,a}}$. We have
 \begin{eqnarray*}
  \sigma_{-+,a} (A)
 & = & \theta \parth f ( -D +a) ( -t)^j
     + \parth \theta f_1 ( -D +a) ( -t)^j  \\
 & = & \theta \parth ( -t)^j f ( -D -j +a)
     + \parth \theta ( -t)^j f_1 ( -D -j +a).
 \end{eqnarray*}
 We obtain by comparing with $ \sigma_{-+,a} (A) = -A$ that
 \begin{eqnarray*}
  - f (D) = ( -1)^j f_1 ( -D -j +a), \quad
  - f_1 (D) = ( -1)^j  f ( -D -j +a).
 \end{eqnarray*}
 These two equations are equivalent.

 On the other hand, we have
 \begin{eqnarray*}
  \sigma_{-+,a} (B)
 & = & g_- ( -D +a) ( -t)^{j -1} \theta - g_+ ( -D +a) ( -t)^j \parth   \\
 & = & \theta ( -t)^{j -1} g_- ( -D -j +1 +a) + \parth ( -t)^j g_+ ( -D -j +a)
 \end{eqnarray*}
 We obtain by comparing with $ \sigma_{-+,a} (B) = -B$ that
 \begin{eqnarray*}
  - g_- (D) & = & ( -1)^{j -1}g_- ( -D -j +1 +a), \\
  g_+ (D) & = & ( -1)^{j} g_+ ( -D -j +a).
 \end{eqnarray*}
 If we let $f (w) = g_- (w + (a +1 -j)/2 )$ and $g(w) = g_+ (w + (a -j)/2 )$,
 we see that $ f (w), g(w) \in \C^{(\overline{j +1})} [w]$.
\end{proof}

In a similar way one determines $\sd^{\sigma_{--,a} }$ as follows.

\begin{prop}\label{zeroform4}
 We have the following description for the graded subspaces
 $\sd_j^{\sigma_{--,a} }$ and $\sd^{\sigma_{--,a} }_{j - 1/2}$ $ (j \in \Z)$
 of \,$\sd^{\sigma_{--,a} }$:
  \begin{eqnarray*}
 \sd^{\sigma_{--,a} }_j
  & = & \left \{
                 \parth \theta t^j f ( D) - ( -1)^j \theta \parth t^j f(-D -j +a)
                 \mid f (w) \in \C [w]
        \right \} ,      \\
  \sd^{\sigma_{--,a} }_{j - 1/2}
  & = & \left \{ \theta t^{j -1} f (D + (j -a -1)/2)
                     + \parth t^j g (D + (j -a)/2) \right. \\
   && \quad\quad\quad\quad\quad\quad\quad\quad\quad\quad\quad\quad\quad
  \left.  \mid  f(w ), g (w)  \in \C^{(\overline{j })} [w]
         \right \}.
 \end{eqnarray*}
\end{prop}

Observe that the $ \sigma_{\pm\pm, a}$'s for different $a$ are
transformed into each other by some spectral flow. More explicitly
we have $ \sigma_{\pm\pm, a} \cdot \Theta_{s} = \sigma_{\pm\pm, a+
s}
   = \Theta_{-s} \cdot \sigma_{\pm\pm, a}$.
Thus the Lie superalgebra $\sd^{\sigma_{\pm \pm, a}}$ is
isomorphic to $\sd^{\sigma_{\pm \pm, a'}}$ for $a, a' \in \C$. As
we shall see, it is most convenient to choose $a = -1$ which we
will fix from now on. Also we will write ${}^{\pm\pm}\sd$ for
$\sd^{\sigma_{\pm\pm,-1}}$. Similarly the corresponding
subalgebras of $\hsd$, induced by the $2$-cocycle
\eqnref{sdcocycle}, will be denoted by $^{\pm\pm}\widehat{\sd}$.

We conclude this section by displaying a distinguished spanning
set for ${}^{\pm\pm}\sd$, that will be of use later on. Let
\begin{eqnarray*}
  {}^{+\pm} W^{n}_{k} & = &  t^k \left(
                    [D]_n \parth \theta - [ -D -k -1]_n \theta\parth
                \right),  \\
  {}^{+\pm} W^{n}_{+,k} & = &  t^{k} \left(
                    [D]_n \mp [ -D -k-1]_n
                \right) \theta, \\
  {}^{+\pm}W^{n}_{ -,k} & = &  t^{k} \left(
                    [D]_n \pm [ -D -k -1]_n
                \right) \parth.
\end{eqnarray*}
Then ${}^{+\pm} W^{n}_{k}, {}^{+\pm} W^{n}_{+,k}$ and
${}^{+\pm}W^{n}_{ -,k}$ $(  k \in \Z, n \in \Z_{+})$ span
${}^{+\pm}\sd$. Similarly we set
\begin{eqnarray*}
  {}^{-\pm} W^{n}_{k} & = & t^k \left(
                    [D]_n \parth \theta + ( -1)^{k+1} [ -D -k -1]_n \theta \parth
                \right),  \\
  {}^{-\pm} W^{n}_{+,k} & = & t^{k} \left(
                    [D]_n \pm ( -1)^{k+1} [ -D -k-1]_n
                \right) \theta, \\
  {}^{-\pm}W^{n}_{ -,k} & = & t^{k} \left(
                    [D]_n \pm ( -1)^{k} [ -D -k -1]_n
                \right) \parth.
\end{eqnarray*}
Then ${}^{-\pm} W^{n}_{k}, {}^{-\pm} W^{n}_{+,k}$ and
${}^{-\pm}W^{n}_{ -,k}$ $( k \in \Z, n \in \Z_{+})$ span
$^{-\pm}\sd$.
\section{Criterion for quasifiniteness of modules over ${}^i\widehat{\sd}$}
\label{sect_criterion}

Given a $\hf \Z$-graded Lie superalgebra $\g = \bigoplus_{j \in
\hf \Z} {\g}_j$ (possibly $\dim {\g}_j = \infty$) with
$\bigoplus_{j \in \Z} {\g}_j$ as the even part and $\bigoplus_{j
\in \hf + \Z} {\g}_j$ as the odd part. A $\mf g$-module $M =
\bigoplus_{j \in \hf \Z} M_j $ is {\em graded} if ${\g}_i M_j
\subset M_{i+j}.$ $M$ is called {\em quasifinite} if $ \dim M_j <
\infty$, following the terminology of Kac and Radul \cite{KR1}. We
assume further that $\g_0$ is abelian and regarded as the Cartan
subalgebra of $\g$.

Given a highest weight $\xi$ in the restricted dual ${\mf g}_0^*$
of $ {\g}_0$, we denote by $L( {\g}; \xi)$ the irreducible
representation of $\g$ with highest weight $\xi$ and a highest
weight vector $v_{\xi}$ annihilated by ${\g}_+ = \bigoplus_{j \in
\hf \mathbb N} {\g}_j$.  That is it is the unique irreducible
quotient of the $\G$-module obtained by inducing from the
one-dimensional $(\bigoplus_{j\ge 0}\G_j)$-module determined by
$\xi$.

From now on $\g$ will always be one of the Lie superalgebras
$\hsdzero$ or ${}^{\pm\pm}\widehat{\sd}$. Note that any of these
Lie superalgebras has infinite-dimensional graded subspaces. Our
problem is to determine for which $\xi$ the corresponding
irreducible highest weight module $L(\G,\xi)$ is quasifinite. Note
that the Verma module is never quasifinite.  We remark that the
$\WW_{1+\infty}$ case was first determined in \cite{KR1} and
subsequently the $\hsd$ case was determined in \cite{AFMO1}. Our
analysis below has a somewhat different flavor.

An obvious necessary condition for $L(\G,\xi)$ to be quasifinite
is that $ {\g}_{ -j} v_{\xi}$ $(j > 0)$ is finite-dimensional,
since ${\g}_{ -j} v_{\xi} \subset L(\G,\xi)_{-j}$. We will be
particularly interested in the case $j = \hf$. We will show
eventually that the condition that $ {\g}_{ -1/2} v_{\xi}$ is
finite dimensional is also sufficient for $L( {\g}; \xi)$ to be
quasifinite.

For $j\in\hf\Z$ we define $P_j=\{X\in\G_j\mid Xv_\xi=0\}$.
Obviously $P=\bigoplus_{j}P_j$ is a parabolic subalgebra of $\G$.

\begin{prop}\label{quasifinite}
Let $\G=\bigoplus_{j\in\Z}\G_j$ be a $\Z$-graded Lie superalgebra
such that $\G_0$ is abelian. Let $\xi\in\G_0^*$ and $L(\G,\xi)$ be
the irreducible highest weight module of highest weight $\xi$ and
$v_\xi$ a corresponding highest weight vector. Suppose that
$\G_jv_\xi$ is finite-dimensional for all $j$.  Then $L(\G,\xi)$
is quasifinite.
\end{prop}

\begin{proof}
Of course if $L(\G,\xi)$ is quasifinite, then $\G_jv_\xi$ is
necessarily finite-dimensional for all $j$.

Conversely suppose that $\G_jv_\xi$ is finite-dimensional for all
$j$.  We need to show that for any positive integers
$i_1,i_2,\cdots,i_m\in\N$ the space
$\G_{-i_1}\G_{-i_2}\cdots\G_{-i_m}v_\xi$ is finite-dimensional. We
will do this by induction on $m$, with $m=1$ being the hypothesis
of the proposition.

Now $\G_{-i_m}v_\xi$ is finite-dimensional, hence there exist
finitely many elements $x_1,x_2,\cdots,x_n\in\G_{-i_m}$ such that
$\G_{-i_m}$ is spanned by $x_1v_\xi,x_2v_\xi,\cdots,x_nv_\xi$ This
implies that $\G_{-i_1}\G_{-i_2}\cdots\G_{-i_m}v_\xi$ lies in the
span of $\G_{-i_1}\G_{-i_2}\cdots\G_{-i_{m-1}}x_iv_\xi$, for $1\le
i\le n$.  Thus it is enough to prove that, for each $i$, each of
these $\G_{-i_1}\G_{-i_2}\cdots\G_{-i_{m-1}}x_iv_\xi$ is
finite-dimensional.

Now we have
\begin{align*}
\G_{-i_1}\G_{-i_2}\cdots\G_{-i_{m-1}}x_iv_\xi\subseteq
x_i\G_{-i_1}\G_{-i_2}\cdots\G_{-i_{m-1}}v_\xi +
[x_i,\G_{-i_1}\G_{-i_2}\cdots\G_{-i_{m-1}}]v_\xi.
\end{align*}
By induction $\G_{-i_1}\G_{-i_2}\cdots\G_{-i_{m-1}}v_\xi$ is a
finite-dimensional vector space and hence
$x_i\G_{-i_1}\G_{-i_2}\cdots\G_{-i_{m-1}}v_\xi$ is
finite-dimensional. Thus it is enough to prove that
$[x_i,\G_{-i_1}\G_{-i_2}\cdots\G_{-i_{m-1}}]v_\xi$ is
finite-dimensional as well.  But
\begin{align*}
[x_i,\G_{-i_1}\G_{-i_2}\cdots\G_{-i_{m-1}}]v_\xi\subseteq&
[x_i,\G_{-i_1}]\G_{-i_2}\cdots
\G_{-i_{m-1}}v_\xi+\G_{-i_1}[x_i,\G_{-i_2}]\cdots\G_{-i_{m-1}}v_\xi\\
&+\cdots\cdots+\G_{-i_1}\G_{-i_2}\cdots[x_i,\G_{-i_{m-1}}]v_\xi.
\end{align*}
However each of these spaces on the right-hand side lies in some
finite-dimensional vector space by induction hypothesis.  Thus
$[x_i,\G_{-i_1}\G_{-i_2}\cdots\G_{-i_{m-1}}]v_\xi$ is also
finite-dimensional.
\end{proof}

\propref{quasifinite} implies that quasifiniteness is equivalent
to the spaces $\G_j/P_j$ being finite-dimensional, for all $j<0$.
Hence our problem at hand reduces to showing that the subspaces
$P_j$ have finite co-dimension in $\G_j$ for all $j<0$.

For a given $\xi$, we introduce the following generating functions
in a variable $x$ for the highest weights:

\begin{eqnarray*}
 {}^{0}\Delta_0 (x)      &=& -\xi (\parth \theta \sinh (xD) ), \\
 {}^0\Delta_1 (x)        &=& -\xi ( \theta\parth \sinh (x( D + 1/2) ),  \\
 {}^{\pm\pm}\Delta_0 (x) &=& -\xi ( \sinh (x (D +1/2)) ), \\
 {}^{\pm\pm}\Delta_1 (x)  &=& -\xi ((\theta \parth -\parth \theta )
 \cosh(x(D+1/2))),
\end{eqnarray*}
where we recall that $\sinh (z) =(e^z -e^{-z})/2$ and $\cosh (z)
=(e^z+e^{-z})/2$.
\subsection{The case ${\g} = \hsdzero$}

We will work out this case in detail. The other cases can be
worked out in a similar and straightforward fashion and we will
omit the arguments for the sake of preserving space.

Let us first recall that for $j\in\Z$
\begin{align*}
 &\G_j={}^0\hsd_j
   = \{ \parth \theta t^j f( D + j/2 ) +
       \theta \parth  t^j g( D + (j +1)/2 )
  \mid f , g \in \C^{\bar{(j)}} [w]  \} ,         \\
  &\G_{j-\hf}={}^0\hsd_{j - 1/2}
   = \{ \theta t^{j -1} g_0 (D) + \parth (-t)^j g_0 ( -D -j),
       g_0 \in \C[w]   \}.
\end{align*}

For $j\in\Z$, we define $I^0_j$ (respectively $I^1_j$) to be the
subspace of $\C^{\bar{(j)}}[w]$ consisting of polynomials $f(w)$
such that $\parth \theta t^j f( D + j/2 )\in P_j$ (respectively
$\theta
\parth  t^j f( D + (j +1)/2 )\in P_j$). Similarly we let $I_{j-\hf}$ be the subspace of
$\C[w]$ consisting of polynomials $f(w)$ such that $\theta t^{j
-1} f(D)+\parth (-t)^j f ( -D -j)\in P_{j-\hf}$.

\begin{lem}\label{lemsdzero}

For $j\in\Z$, $I_{j-\hf}$ is an ideal of $\C[w]$, while $I^0_j$
and $I^1_j$ are $\C^{(0)}[w]$-submodules of $\C^{\bar{(j)}}[w]$.
Thus there exists monomials $\alpha_{j-\hf}(w)\in\C[w]$ and
$\beta^0_j(w),$ $ \beta^1_j(w)\in\C^{\bar{(j)}}[w]$ such that
$I_{j-\hf}=\C[w]\alpha_{j-\hf}(w)$,
$I^0_j=\C^{(0)}[w]\beta^0_j(w)$ and
$I^1_j=\C^{(0)}[w]\beta^1_j(w)$.  Furthermore $$P_j = \{ \parth
\theta t^j f(D+\frac{j}{2}) +
       \theta \parth  t^j g( D + \frac{j+1}{2})
  \mid f\in I^0_j , g \in I^1_j\}.$$
\end{lem}

\begin{proof}
Consider $f(w)\in\C^{(1)}[w]$ so that $\parth\theta f(D)\in\G_0$.
Take any $\theta t^{j-1}g(D)+\parth (-t)^jg(-D-j)\in P_{j-\hf}$,
for some $g(w)\in I_{j-\hf}$.  We compute their bracket
\begin{align*}
[\parth\theta f(D),\; &\theta t^{j-1}g(D)+\parth (-t)^jg(-D-j)]\\
=&-(\theta t^{j-1}g(D)f(D)+\parth(-t)^j f(-D-j)g(-D-j))\in
P_{j-\hf}.
\end{align*}
Thus $f(w)g(w)\in I_{j-\hf}$, for any $f(w)\in\C^{(1)}[w]$, that
is, $I_{j-\hf}$ is invariant under the multiplication of any odd
polynomial. But then it follows that $I_{j-\hf}$ is invariant
under the multiplication of any polynomial, and hence is an ideal
of $\C[w]$.

Again take any $f(w)\in\C^{(1)}[w]$ so that $\parth\theta f(D) \in
\G_0$, and any element $$\parth\theta
t^jg(D+\frac{j}{2})+\theta\parth t^j\bar{g}(D+\frac{j+1}{2})\in
P_{j}$$ for some $g(w) \in I^0_j$ and $\bar{g}(w) \in I^1_j$. We
compute

\begin{align*}
[\parth\theta f(D),\parth\theta t^jg(D+\frac{j}{2})+\theta\parth
t^j\bar{g}(D+\frac{j+1}{2})]=
\parth\theta g(D+\frac{j}{2})h(D+\frac{j}{2}),
\end{align*}
where $h(D+\frac{j}{2})=f(D+j)-f(D)$.  As $f(w)$ ranges over all
odd polynomials, $h(w)$ ranges over all even polynomials.  From
this it follows that $I^0_j$ is a module over $\C^{(0)}[w]$.

Replacing $\parth\theta f(D)$ with $\theta\parth f(D)$ in the
above calculation one shows similarly that $I^1_j$ is a module
over $\C^{(0)}[w]$.
\end{proof}

It follows from \lemref{lemsdzero} that that $L(\G,\xi)$ is
quasifinite if and only if $I^0_j, I^1_j$, and $I_{j+\hf}$ are all
nonzero for $j\in -\N$.

\begin{lem} Suppose that $I_{-\hf}\not=0$ and $j\in \N$. We have:
\begin{itemize}
\item[(i)] $I^a_j\not=0$ if and only if $I_{j-\hf}\not=0$, for $a=0,1$.
\item[(ii)] $I_{j+\hf}\not=0$ if and only if $I^a_{j}\not=0$, for $a=0,1$.
\end{itemize}
\end{lem}

\begin{proof}  We will first prove (ii).
Suppose that $0\not=f(w)\in I^0_j$. Let $g(w)\in\C[w]$ so that
$\theta g(D)-\parth t g(-D-1)\in\G_{\hf}$.  We compute
\begin{align*}
[\parth\theta t^jf(D+\frac{j}{2})&,\theta g(D)-\parth t
g(-D-1)]=\\ &\theta t^j h(D)+\parth (-t)^{j+1}h(-D-j-1)\in
P_{j+\hf},
\end{align*}
where $h(D)=g(D+j)f(D+\frac{j}{2})$.  Since $f(w)\not=0$ and
$g(w)$ can be taken to be an arbitrary non-zero element of
$\C[w]$, it follows that $I_{j+\hf}\not=0$. The claim that
$I^1_j\not=0$ implies that $I_{j+\hf}\not=0$ is proved similarly
and thus omitted.

Conversely suppose that $I_{j+\hf}\not=0$.  Let $0\not=f(w)\in
I_{j+\hf}$ and $0\not=g(w)\in I_{-\hf}$ and compute
\begin{align*}
[\theta t^{-1}g(D)+\parth g(-D),\theta t^j
f(D)+\parth(-t)^{j+1}&f(-D-j-1)]=\\ &\parth\theta
t^jh(D+\frac{j}{2}) + \theta\parth t^j k(D+\frac{j+1}{2}),
\end{align*}
where $h(D+\frac{j}{2})=f(D)g(-D-j)+(-1)^{j+1}f(-D-j)g(D)$ and
$k(D+\frac{j+1}{2})=f(D)g(-D) + (-1)^{j+1} f(-D-j-1)g(D+j+1)$.

Now $h(w)=0$ if and only if
$f(w-\frac{j}{2})g(-w-\frac{j}{2})\in\C^{\overline{(j+1)}}[w]$.
Now since $I_{j+\hf}$ is an ideal of $\C[w]$, we may replace
$f(w)\in I_{j+\hf}$ above by $0\not=(w+\frac{j}{2})f(w)\in
I_{j+\hf}$. This implies then that
$wf(w-\frac{j}{2})g(-w-\frac{j}{2})\in\C^{\overline{(j+1)}}[w]$ as
well. But this is impossible unless
$f(w-\frac{j}{2})g(-w-\frac{j}{2})=0$, which is a contradiction.
Hence $I^0_{j}\not=0$. Similarly $I^1_j\not=0$.  This proves (ii).

Next we prove (i). For this suppose $0\not=f(w)\in I^0_{j}$ and
$0\not=g(w)\in I_{-\hf}$. Then
\begin{align*}
[\parth\theta t^jf(D+\frac{j}{2}),\theta t^{-1}g(D)+\parth
g(-D)]=\\ (-1)^j\Bigl{(}\theta
t^{j-1}g(D+j)f(-D-\frac{j}{2})+(-t)^j\parth
f(D+\frac{j}{2})g(-D)\Bigl{)}\in P_{j-\hf},
\end{align*}
and hence $I_{j-\hf}\not=0$.  Similarly $I^1_j\not=0$ implies
$I_{j-\hf}\not=0$.

Finally suppose that $0\not=g(w)\in I_{\hf}=\C[w]$ and
$0\not=f(w)\in I_{j-\hf}$ and consider

\begin{align*}
[\theta g(D)-\parth t g(-D-1),&\theta
t^{j-1}f(D)+\parth(-t)^jf(-D-j)]=\\ &\parth\theta
t^jh(D+\frac{j}{2})+ \theta\parth t^jk(D+\frac{j+1}{2}),
\end{align*}
where $h(D+\frac{j}{2})=(-1)^jf(-D-j)g(D)-f(D)g(-D-j)$ and
$k(D+\frac{j+1}{2})=(-1)^jg(D+j)f(-D-j)-f(D+1)g(-D-1)$.  Now the
same argument proving $I_{j+\hf}\not=0$ implying $I^a_j\not=0$ can
be used to derive that $I^a_j\not=0$, for $a=0,1$.
\end{proof}

\begin{cor}\label{sdzerocor1}
The irreducible highest weight module $L(\G,\xi)$ is quasifinite
if and only if $I_{-\hf}\not=0$.
\end{cor}

In particular setting $\alpha(w)=\alpha_{-\hf}(w)$ in
\lemref{lemsdzero} we obtain the following.

\begin{cor}\label{sdzerocor2}
 There exists a monomial $\alpha(w) \in \C[w]$ such that
 any element in $(\hsdzero)_{-\hf}$ annihilating $v_{\xi}$ is of the form
 $\theta t^{-1} g(D)\alpha(D) + \parth g( -D)\alpha( -D)$, for
 some $g(w)\in\C[w]$.
\end{cor}

Take this $\theta t^{-1} \alpha (D) + \parth \alpha ( -D) \in {\mf
g}_{- 1/2}$ and any element $ \theta f(D) - \parth t f( -D -1) \in
{\g}_{1/2}$. One computes
\begin{eqnarray*}
      \begin{array}{rcl}
      \lefteqn{
            [\theta f(D) - \parth t f( -D -1),
                 \theta t^{-1} \alpha (D) +  \parth \alpha ( -D)] } \\
     & = & \parth \theta ( \alpha ( -D) f(D) - \alpha (D) f( -D) )  \\
     & & + \theta \parth (\alpha ( -D) f(D) - \alpha (D +1) f( -D
     -1))
    - \alpha (0) f(0) C.
      \end{array}
\end{eqnarray*}
Since $[\theta f(D)-\parth tf(-D-1),\theta t^{-1}\alpha(D)+\parth
\alpha(-D)]v_\xi=0$, we obtain
\begin{eqnarray*}
      \begin{array}{rcl}
      \lefteqn{
             \xi \left\{ \parth \theta ( \alpha ( -D) f(D) - \alpha (D) f( -D) )
      \right. } \\
     && \left.
    + \theta \parth (\alpha ( -D) f(D) - \alpha (D +1) f( -D -1)
    - \alpha (0) f(0) c \right\} = 0.
      \end{array}
\end{eqnarray*}

Putting $f(D) = e^{- xD}$ this is equivalent to
\begin{equation}\label{eqnsdzero}
\alpha \left(\frac{d}{dx} \right) \left( {}^0\Delta_0 (x)
  +   {}^0\Delta_1 (x)e^{x/2}
  -  \frac{c}{2} \right) = 0.
\end{equation}

Combining the results of this section we obtain the following
characterization of quasifinite $\hsdzero$-modules.

\begin{thm} The irreducible highest weight module
$L(\hsdzero,\xi)$ is quasifinite if and only if ${}^0\Delta_0(x)$
and ${}^0\Delta_1(x)$ satisfy the differential equation
\eqnref{eqnsdzero}.
\end{thm}
\subsection{The case ${\g} = {}^{\pm\pm}\hsd$}

First we let ${\g} = \hsdpp$. The arguments leading to the
following proposition is similar to the ones we have given in the
previous section leading to \corref{sdzerocor1} and
\corref{sdzerocor2} and so we will omit the details.

\begin{prop}
  There exist $\alpha (w) \in \C^{(1)} [w]$ and $\beta(w) \in \C^{(0)} [w]$
 such that any element in $(\hsdpp)_{-\hf}$ annihilating $v_{\xi}$ is of the form
 $\theta t^{-1} f(D)\alpha(D) + \parth g( D)\beta(D+1/2),$ where
 $f(w), g (w) \in \C^{(0)} [w]$.  Furthermore the irreducible
 highest weight module $L(\G,\xi)$ is quasifinite if and only if
 $\alpha(w)\not=0$ and $\beta(w)\not=0$.
\end{prop}

 Take $\theta f( D + 1/2) + \parth t g( D +1) \in {\g}_{ 1/2}$
and $ \theta t^{ -1} \alpha (D) + \parth \beta (D + 1/2) \in P_{
-1/2}$, where $\alpha, f \in \C^{(1)} [w]$ and $\beta, g \in
\C^{(0)} [w]$. One computes

\begin{eqnarray*}
  \begin{array}{rcl}
   \lefteqn{
    [ \theta f( D + 1/2) + \parth t g( D +1),
   \theta t^{ -1} \alpha (D) + \parth \beta (D + 1/2)] } \\
    & = & \beta (D + 1/2) f (D + 1/2)
    + \hf ( \alpha (D +1) g( D +1) + \alpha (D) g(D) )   \\
    & & + \hf (\theta \parth - \parth \theta )
              ( \alpha (D +1) g( D +1) - \alpha (D) g(D) )
     + \alpha (0) g (0) C.
  \end{array}
\end{eqnarray*}

This implies that
\begin{align}\label{pdeqn0}
  \xi&\Bigl{(}
  \beta (D + 1/2) f (D + 1/2)
    + \hf ( \alpha (D +1) g( D +1) + \alpha (D) g(D) ) \\
   &+\hf (\theta \parth - \parth \theta )
              ( \alpha (D +1) g( D +1) - \alpha (D) g(D) )
     + \alpha (0) g (0) c \Bigl{)} = 0.\nonumber
\end{align}
By putting $g(D)=0$ and $f(D+1/2)=\hf(e^{x(D +
1/2)}-e^{-x(D+1/2)})$ we have

\begin{align}\label{pdeqn1}
\beta(\frac{d}{dx})\;\;({}^{++}\Delta_0(x) )=0.
\end{align}
Now putting $f(D) = 0$ and $g(D) = \hf (e^{xD} + e^{ -xD})$ we
obtain

\begin{align}\label{pdeqn2}
\alpha (\frac{d}{dx}) \left( (e^{x/2} + e^{ -x/2} )\cdot
{}^{++}\Delta_0 (x) +
 (e^{x/2} - e^{ -x/2} )\cdot{}^{++}\Delta_1 (x) - 2c \right)
  = 0.
\end{align}

\vspace{.3cm}

 Next, we consider ${\g} = \hsdpm$. In this case we
have similarly the following proposition.

\begin{prop}
  There exist  $\alpha (w) \in \C^{(0)} [w]$ and
 $\beta(w) \in \C^{(1)} [w]$ such that
 any element in $(\hsdpm)_{-\hf}$ annihilating $v_{\xi}$ is of the form
 $\theta t^{-1} f(D)\alpha(D) + \parth g( D)\beta(D+1/2),$ where
 $f(w), g (w) \in \C^{(0)} [w]$. Furthermore the irreducible
 highest weight module $L(\G,\xi)$ is quasifinite if and only if
 $\alpha(w)\not=0$ and $\beta(w)\not=0$.
\end{prop}

 Take $\theta f( D + 1/2) + \parth t g( D +1) \in {\g}_{ 1/2}$
and $ \theta t^{ -1} \alpha (D) + \parth \beta (D + 1/2) \in
{\g}_{ -1/2}$, where $\alpha, f \in \C^{(0)} [w]$ and $\beta, g
\in \C^{(1)} [w]$. As before we obtain identity \eqnref{pdeqn0}.
Now using the pairs $g(D)=0$ together with $f(D+1/2)=\hf(e^{x(D +
1/2)}+e^{-x(D+1/2)})$ and $f(D) = 0$ together with $g(D) = \hf
(e^{xD} + e^{ -xD})$ we obtain again the equations \eqnref{pdeqn1}
and \eqnref{pdeqn2} where $\Delta_i^{++}(x)$ is replaced by
$\Delta_i^{+-}(x)$.

\vspace{.3cm}

Next, we consider ${\g} = \hsdmp$. In this case we have similarly
the following proposition.

\begin{prop}
  There exist $\alpha(w),\beta (w) \in \C^{(0)} [w]$ such that
 any element in $(\hsdpm)_{-\hf}$ annihilating $v_{\xi}$ is of the form
 $\theta t^{-1} f(D)\alpha(D) + \parth g( D)\beta(D+1/2),$ where
 $f(w), g (w) \in \C^{(0)} [w]$. Furthermore the irreducible
 highest weight module $L(\G,\xi)$ is quasifinite if and only if
 $\alpha(w)\not=0$ and $\beta(w)\not=0$.
\end{prop}

 Take $\theta f( D + 1/2) + \parth t g( D +1) \in {\g}_{ 1/2}$
and $ \theta t^{ -1} \alpha (D) + \parth \beta (D + 1/2) \in
{\g}_{ -1/2}$, where $\alpha, f \in \C^{(0)} [w]$ and $\beta, g
\in \C^{(1)} [w]$, and taking their bracket gives identity
\eqnref{pdeqn0}.
Now using the pairs $g(D)=0$ together with $f(D+1/2)=\hf(e^{x(D +
1/2)}-e^{-x(D+1/2)})$ and $f(D) = 0$ together with $g(D) = \hf
(e^{xD} - e^{ -xD})$ we obtain again \eqnref{pdeqn1} and
\eqnref{pdeqn2} where $\Delta_i^{++}(x)$ is replaced by
$\Delta_i^{-+}(x)$.

\vspace{.3cm}

Finally, we consider ${\g} = \hsdmm$. In this case we have
similarly the following proposition.

\begin{prop}
  There exist $\alpha (w), \beta(w) \in \C^{(1)} [w]$ such that
 any element in $(\hsdpm)_{-\hf}$ annihilating $v_{\xi}$ is of the form
 $\theta t^{-1} f(D)\alpha(D) + \parth g( D)\beta(D+1/2),$ where
 $f(w), g (w) \in \C^{(0)} [w]$. Furthermore the irreducible
 highest weight module $L(\G,\xi)$ is quasifinite if and only if
 $\alpha(w)\not=0$ and $\beta(w)\not=0$.
\end{prop}

As before we obtain identity \eqnref{pdeqn0}. Now we use the pairs
$g(D)=0$ together with $f(D+1/2)=\hf(e^{x(D +
1/2)}+e^{-x(D+1/2)})$ and $f(D) = 0$ together with $g(D) = \hf
(e^{xD} + e^{ -xD})$ to obtain \eqnref{pdeqn1} and
\eqnref{pdeqn2}, where $\Delta_i^{++}(x)$ is replaced by
$\Delta_i^{--}(x)$.

We summarize our results in the following.

\begin{thm} The irreducible highest weight module
$L(^{\pm\pm}\widehat{\sd},\xi)$ is quasifinite if and only if
${}^{\pm\pm}\Delta_0(x)$ and ${}^{\pm\pm}\Delta_1(x)$ satisfy the
following differential equations:

 \begin{eqnarray*}
  \beta(\frac{d}{dx})({}^{\pm\pm}\Delta_0(x)) &=& 0, \\
  \alpha (\frac{d}{dx}) \left( (e^{x/2} + e^{ -x/2})\cdot {}^{\pm\pm}\Delta_0 (x) \right)
   +  \alpha (\frac{d}{dx}) \left( (e^{x/2} - e^{ -x/2} ) \cdot
   {}^{\pm\pm}\Delta_1 (x) - 2c \right)
   &=& 0.
 \end{eqnarray*}
\end{thm}
\section{Subalgebras of Lie superalgebras $\hgltwo$ and
$\hsd$}\label{sect_embed}

In this section we will describe Lie subalgebras of $\gltwo$ of
types $OSP$ and $P$. We study the embedding of the Lie subalgebras
of $\sd$ introduced in Section~\ref{sect_subalg} into (the
subalgebras of) $\gltwo$. We then study the relations among the
central extensions of these superalgebras. Also cf. \cite{AFMO1}
for the embedding of $\hsd$ into $\hgltwo$.
\subsection{Lie superalgebra $\gltwo$}

Let $M_{\infty}$ be the associative algebra consisting of matrices
$(a_{ij})_{i,j \in \Z}$ such that $a_{ij} = 0$ for $ |i -j | >>
0$. We denote by $\glone$ the Lie algebra obtained from
$M_{\infty}$ by taking the usual commutator. Denote by $E_{ij}$
the elementary matrix with $1$ at the $(i,j)$-th entry and $0$
elsewhere. Define the Lie superalgebra $\gltwo$ to be $\glone
\otimes M (1,1)$ with the induced $\Z_2$-graded structure from $M
(1,1)$ (recall that $M(1,1)$ is defined earlier in
Section~\ref{sect_classical}).

One may have two different ways of looking at $\gltwo$. First we
may regard
 $\gltwo = \bigoplus_{a = 0, 1, \pm} \glone M_a $, that is,
 \begin{eqnarray*}
 \gltwo =
 \left[
  \begin{array}{cc}
   \glone   & \glone       \\
   \glone   & \glone
 \end{array} \right].
   \label{mat_inf}
 \end{eqnarray*}
Secondly, one may identify $\gltwo = \{ ( a_{ij})_{i,j \in \Z/2}
\mid a_{ij} = 0,\ {\rm for}\ | i -j | >>0 \}$. Under this
identification, the $\Z_2$-graded structure is given by letting
$E_{ij}$ $(i -j \in \Z)$ be even while $E_{ij}$ $(i -j \in \Z +
1/2)$ be odd. One could easily identify the two presentations of
$\gltwo$ as follows ($i, j \in \Z$):
\begin{eqnarray*}
 E_{ij} M_0 = E_{ij}, & E_{ij} M_1 = E_{i- 1/2, j- 1/2}, \\
 E_{ij} M_+ = E_{i,j - 1/2}, & E_{ij} M_-= E_{i- 1/2, j}.
\end{eqnarray*}

The Lie superalgebra $\gltwo$ is equipped with a natural
$\hf\Z$-gradation
 $$\gltwo=\bigoplus_{r\in\hf\Z}(\gltwo)_{r}  $$
where $(\gltwo)_{r}$ is the completion of the linear span of
$E_{ij}$ with $j-i=r$ in the notation of the second realization
above. We will call this gradation the {\em principal gradation}
of $\gltwo$. This gradation naturally gives a choice of a Borel
subalgebra $\btwo=\bigoplus_{r\ge 0}(\gltwo)_{r}$ and a choice of
a Cartan subalgebra $\htwo= (\gltwo)_0$.  We obtain thus a
triangular decomposition of $\gltwo$
\begin{equation*}
\gltwo=(\gltwo)_+\bigoplus(\gltwo)_0\bigoplus(\gltwo)_-,
\end{equation*}
where $(\gltwo)_\pm=\bigoplus_{r\in\hf\N}(\gltwo)_{\pm r}$.
\subsection{Subalgebras $\ospd$ and $\pd_{\pm\pm}$ of $\gltwo$}

Consider the superspace $\C^{\infty|\infty}$ with even basis
vectors $\{e_i|i\in\Z\}$ and odd basis vectors
$\{e_r|r\in\hf+\Z\}$. The Lie superalgebra $\gltwo$ acts on
$\C^{\infty|\infty}$ in the usual way:
\begin{equation*}
E_{ij}e_k=\delta_{jk}e_i,\quad i,j,k\in\hf\Z.
\end{equation*}

Consider the {\em even super-symmetric} non-degenerate bilinear
form $(\cdot|\cdot)$ on $\C^{\infty|\infty}$ given by
\begin{eqnarray*}
(e_i|e_j)&=&(-1)^i\delta_{i,-j},\quad i,j\in\Z;\\
(e_r,e_s)&=&(-1)^{r+\hf}\delta_{r,-s},\quad r,s\in\hf+\Z; \\
(e_i|e_r)&=&0,\quad i\in\Z, r\in\hf+\Z.\\
\end{eqnarray*}
Associated to this form we define the Lie superalgebra $\ospd
=\ospd_0 \oplus \ospd_1$ to be the subalgebra of $\gltwo$
preserving this form, where
$$\ospd_\epsilon=\{T\in(\gltwo)_{\epsilon}|(Tv|w)=-(-1)^{\epsilon|v|}(v,Tw)\},
 \;\;\epsilon =0,1.
$$
This is a Lie superalgebra of type $OSP$. It is easy to see that
the subalgebra $\ospd^f$ spanned by the following elements is a
dense subalgebra inside $\ospd$ ($i,j\in\Z$, $r,s\in\hf+\Z$):
\begin{equation*}
E_{i,j}-(-1)^{i+j}E_{-j,-i}, \;\; E_{i,r}-(-1)^{i+r-\hf}E_{-r,-i},
\;\; E_{r,s}-(-1)^{r+s-1}E_{-s,-r}.
\end{equation*}
In our first realization of $\gltwo$ the subalgebra $\ospd^f$ is
the span of elements of the form
$(E_{i, j}  - ( -1)^{i+j} E_{ -j, -i} )M_0,
   (E_{i, j} - ( -1)^{i+j} E_{1 -j, 1 -i} )M_1,
$
and $E_{i, j} M_+ + (-1)^{i +j} E_{ 1-j, -i} M_-.$

In the remainder of this subsection we will introduce four
subalgebras of $\hsd$ of type $P$. Note that the type $P$ Lie
superalgebra is one of the two {\em strange} series of Lie
superalgebras, cf.~e.g.~\cite{K} for the finite-dimensional case.

First, we introduce an {\em odd symmetric} non-degenerate bilinear
form $(\cdot|\cdot)$ on $\C^{\infty|\infty}$ which is uniquely
determined by

\begin{equation}\label{oddform1}
(e_i,e_{j-\hf})=\delta_{i+j-\hf,\hf},\quad i,j\in\Z.
\end{equation}
We denote by $\pd_{++}$ the subalgebra of $\gltwo$ preserving this
form.  Then using the notation in the second realization of
$\gltwo$ the subalgebra spanned by

\begin{equation*}
E_{ij}-E_{-j+\hf,-i+\hf},\;\; E_{i,j-\hf}+E_{-j+1,-i+\hf}, \;\;
E_{i-\hf,j}-E_{-j+\hf,-i+1}
\end{equation*}
forms a dense subalgebra $\pd_{++}^{f}$ inside $\pd_{++}$. In the
first realization of $\gltwo$ the subalgebra $\pd_{++}^f$ is the
span of elements of the form
\begin{eqnarray*}
E_{ij} M_1 - E_{1- j, 1- i} M_0,\;\;
    ( E_{ij} - E_{1- j, 1- i}) M_-, \;\;
    ( E_{ij} + E_{1- j, 1- i}) M_+.
\end{eqnarray*}

We may also use \eqnref{oddform1} to (uniquely) determine an {\em
odd skew-symmetric} non-degenerate bilinear form $(\cdot|\cdot)$
on $\C^{\infty|\infty}$. We denote by $\pd_{+-}$ the subalgebra of
$\gltwo$ preserving this form. Then the subalgebra spanned by

\begin{equation*}
E_{ij}-E_{-j+\hf,-i+\hf},\;\; E_{i,j-\hf}-E_{-j+1,-i+\hf}, \;\;
E_{i-\hf,j}+E_{-j+\hf,-i+1}
\end{equation*}
forms a dense subalgebra $\pd_{+-}^{f}$ inside $\pd_{+-}$. In the
first realization $\pd_{+-}^{f}$ corresponds to the subalgebra
spanned by elements of the form
\begin{eqnarray*}
E_{ij} M_1 - E_{1- j, 1- i} M_0, \;\;
    ( E_{ij} + E_{1- j, 1- i}) M_-,  \;\;
    ( E_{ij} - E_{1- j, 1- i}) M_+.
\end{eqnarray*}

Next, we consider the {\em odd symmetric} non-degenerate bilinear
form $(\cdot|\cdot)$ on $\C^{\infty|\infty}$ uniquely determined
by
\begin{equation}\label{oddform2}
(e_i,e_{j-\hf})=(-1)^i\delta_{i+j-\hf,\hf},\quad i,j\in\Z.
\end{equation}
Let $\pd_{-+}$ be the subalgebra of $\gltwo$ preserving this form.
Then the subalgebra spanned by
$E_{ij}-(-1)^{i+j}E_{-j+\hf,-i+\hf},
E_{i,j-\hf}+(-1)^{i+j}E_{-j+1,-i+\hf}, $
and $E_{i-\hf,j}-(-1)^{i+j}E_{-j+\hf,-i+1} $ forms a dense
subalgebra $\pd_{-+}^{f}$ inside $\pd_{-+}$. In the first
realization of $\gltwo$ the subalgebra $\pd_{-+}^{f}$ is the span
of elements of the form
\begin{eqnarray*}
E_{ij} M_1 - (-1)^{i+j}E_{1- j, 1- i} M_0, \quad ( E_{ij} -
(-1)^{i+j}E_{1- j, 1- i}) M_-,
\end{eqnarray*}
and $ ( E_{ij} + (-1)^{i+j}E_{1- j, 1- i}) M_+.$

Similarly we may consider the {\em odd skew-symmetric}
non-degenerate bilinear form $(\cdot|\cdot)$ on
$\C^{\infty|\infty}$ determined by \eqnref{oddform2}.  Let
$\pd_{--}$ be the subalgebra of $\gltwo$ preserving this form.
Then the subalgebra spanned by
\begin{equation*}
E_{ij}-(-1)^{i+j}E_{-j+\hf,-i+\hf}, \quad
E_{i,j-\hf}-(-1)^{i+j}E_{-j+1,-i+\hf},
\end{equation*}
and $E_{i-\hf,j}+(-1)^{i+j}E_{-j+\hf,-i+1}$ forms a dense
subalgebra $\pd_{--}^{f}$ inside $\pd_{--}$. In the first
realization $\pd_{--}^{f}$ corresponds to the subalgebra spanned
by elements of the form
\begin{eqnarray*}
E_{ij} M_1 - (-1)^{i+j}E_{1- j, 1- i} M_0, \quad
    ( E_{ij} + (-1)^{i+j}E_{1- j, 1- i}) M_-,
\end{eqnarray*}
and $( E_{ij} - (-1)^{i+j}E_{1- j, 1- i}) M_+. $ We note the Lie
superalgebras $\ospd$ and $\pd_{\pm\pm}$ inherit from $\gltwo$ the
principal gradation and the triangular decomposition.
\subsection{Embedding of subalgebras of $\sd$ into $\gltwo$}

Take a basis in $\C^{ \infty| \infty} = \C^{ \infty} \bigoplus
\C^{ \infty} \theta
 = t^{s} \C [t, t^{-1}] \bigoplus \theta t^{s} \C [t, t^{-1}] $
as $\{ v_i = t^{-i + s}, v_{ i - 1/2} = t^{ -i + s} \theta, i \in
\Z \}$, where $s\in \C$. The Lie superalgebra $\gltwo$ acts on
$\C^{ \infty| \infty}$ by letting $E_{ij} v_k = \delta_{jk} v_i$,
$i,j,k \in \Z/2$. The Lie superalgebra $\sd$ acts on $\C^{ \infty|
\infty}$ as differential operators. In this way we obtain a family
of embeddings $\phi_{s}$ of $\sd $ into $\gltwo$:
 \begin{eqnarray} \label{eqembed}
 \phi_{s} \left( t^k F(D) \right) =
 \left[
  \begin{array}{cc}
   \phi_{s} \left( t^k f_0 (D)\right)   & \phi_{s} \left( t^k f_+ (D)\right)     \\
   \phi_{s} \left( t^k f_- (D) \right)  & \phi_{s} \left( t^k f_1 (D)\right)
 \end{array} \right],
 \end{eqnarray}
where
\begin{eqnarray*}
 \phi_{s} (t^k f_0 (D))
            & = & \sum_{j \in \Z} f_0 ( -j +s) E_{j -k, j} M_0, \\
 \phi_{s} (t^k f_1 (D))
            & = & \sum_{j \in \Z} f_1 ( -j +s) E_{j -k, j} M_1, \\
 \phi_{s} (t^k f_+ (D))
            & = & \sum_{j \in \Z} f_+ ( -j +s) E_{j -k, j} M_+, \\
 \phi_{s} (t^k f_- (D))
            & = & \sum_{j \in \Z} f_- ( -j +s) E_{j -k, j} M_-.
\end{eqnarray*}
Note that the principal gradation on $\gltwo$ is compatible with
that on $\sd$ under the map $\phi_{s}$.

Denote by $\co $ the algebra of all holomorphic functions on $\C$
with topology of uniform convergence on compact sets. For an
integer $k$ we let
\begin{eqnarray*}
  \co^{(2k)} \equiv \co^{(0)}& = & \{ f \in {\co} \mid f(w) = f(-w) \}, \\
  \co^{(2k+1)} \equiv \co^{(1)} & = & \{ f \in {\co} \mid f(w) = - f( -w) \}.
\end{eqnarray*}
We define a completion $\sd^{\co}$ of $\sd$ consisting  of all
differential operators of the form $ \sum_{a, j} t^j f_a(D) M_a$
where $f_a \in \co$ and $j \in \Z$. This induces completions
$\sdzero^{\co}$ and ${}^{\pm\pm}\sd^{\co}$ of $\sdzero$ and
${}^{\pm\pm}\sd$, respectively. The embedding $\phi_s$ extends
naturally to $\sdzero^{\co}$ and ${}^{\pm\pm}\sd^{\co}$.

We define
\begin{equation*}
I_s:=\{f\in\co|f(n+s)=0,\forall n\in\Z\},
\end{equation*}
and, for $j,k\in\Z$,
\begin{equation*}
 I^{(j)}_{s,k}  =  \{ f \in {\co}^{(j)}
    \mid f (n + k/2 + s) =0,  \forall n \in \Z\},
\end{equation*}

Setting
$J_s=\bigoplus_{k\in\Z}\{\sum_{a=0,1,\pm}t^kf_a(D)M_a|f_a\in
I_s\}$ we have the following proposition whose proof is a
straightforward calculation (compare \cite{KR1}).

\begin{prop} \label{propexactseq}
We have the following exact sequence of Lie superalgebras:

 $$ 0 \longrightarrow J_{ s}\longrightarrow \sd^{\co}
     \stackrel{\phi_s}{\longrightarrow} \gltwo \longrightarrow 0.
 $$
\end{prop}

For $k\in\Z$, let ${}^0J_{s,k}$ be the linear span of
\begin{eqnarray*}
t^k f_0(D + \frac{k}{2})M_0,\;\; t^kf_1(D+\frac{k+1}{2})M_1,\;\;
(-t)^kg(-D-k)M_+ +t^{k-1}g(D)M_-,
\end{eqnarray*}
where $f_0\in I^{(k+1)}_{s,k}$, $f_1\in I^{(k+1)}_{s,k+1}$ and
$g\in I_s$. Set ${}^0J_s=\bigoplus_{k\in\Z}{}^0J_{s,k}$ we can
easily show the following.

\begin{prop} \label{propexactseq2}
  We have the following exact sequences of Lie superalgebras:

\begin{align*}
&0 \longrightarrow {}^0J_{ s}\longrightarrow \sdzero^{\co}
     \stackrel{\phi_s}{\longrightarrow} \gltwo\longrightarrow
     0,\quad s\not\in\hf\Z,\\
&0 \longrightarrow {}^0J_0\longrightarrow \sdzero^\co
     \stackrel{\phi_0}{\longrightarrow}\ospd\longrightarrow 0.
\end{align*}
\end{prop}

For $k\in\Z$ we let ${}^{+\pm}J_{s,k}$ be the linear span of
%
$t^k (f(D+\frac{k}{2})M_0 -f(-D-k-1)M_1),
t^{k-1}f_{-}(D+\frac{k}{2})M_{-}, $ and $
t^kf_+(D+\frac{k-1}{2})M_+,$
where $f\in I_s$ and $f_-\in I^{(1)}_{s,k}$, $f_+\in
I^{(0)}_{s,k-1}$ in the case of ${}^{++}J_{s,k}$, and $f_-\in
I^{(0)}_{s,k}$, $f_+\in I^{(1)}_{s,k-1}$ in the case of
${}^{+-}J_{s,k}$.

For $k\in\Z$ we let ${}^{-\pm}J_{s,k}$ be the linear span of
\begin{eqnarray*}
t^k (f(D+k/{2})M_0 -(-1)^kf(-D-k-1)M_1), \quad t^{k-1}f_{-}(D+
k/{2})M_{-},
\end{eqnarray*}
and $t^kf_+(D+(k-1)/{2})M_+,$ where $f\in I_s$ and $f_-\in
I^{(k)}_{s,k}$, $f_+\in I^{(k)}_{s,k-1}$ in the case of
${}^{-+}J_{s,k}$, and $f_-\in I^{(k+1)}_{s,k}$, $f_+\in
I^{(k+1)}_{s,k-1}$ in the case of ${}^{--}J_{s,k}$.

Similarly, putting
${}^{\pm\pm}J_s=\bigoplus_{k\in\Z}{}^{\pm\pm}J_{s,k}$, we can show
the following.

\begin{prop} \label{propexactseq3}
  We have the following short exact sequences of Lie superalgebras:
\begin{align*}
&0 \longrightarrow {}^{\pm\pm}J_{s} \longrightarrow
 {}^{\pm\pm}\sd^\co
     \stackrel{\phi_s}{\longrightarrow} \gltwo \longrightarrow 0, \quad s \not\in
     \hf\Z,\\
&0 \longrightarrow {}^{\pm\pm}J_{0} \longrightarrow
 {}^{\pm\pm}\sd^\co
     \stackrel{\phi_0}{\longrightarrow} \pd_{\pm \pm} \longrightarrow 0.
     \end{align*}
\end{prop}

\begin{rem}
One can also explicitly describe the images of $\phi_s$ ($s\in
\hf\Z$) which are Lie subalgebras of $\hgltwo$. For our purpose in
this paper, we only need the surjectivity of the homomorphisms
$\phi_s$ in \propref{propexactseq} and \propref{propexactseq2}.
\end{rem}
\subsection{Lifting of the embeddings to the central extensions}

Next we extend the above homomorphism $\phi_{s}$ to a homomorphism
between the central extensions of the corresponding Lie
superalgebras.

We begin by describing the central extension of $\gltwo$. The Lie
superalgebra $\hgltwo$ is the central extension of $\gltwo$ by a
1-dimensional vector space spanned by ${\bf 1}$. The $2$-cocycle
giving rise to this central extension is explicitly given by

\begin{equation*}
\alpha(A,B)={\rm Str}([J,A]B),\quad A,B\in\gltwo,
\end{equation*}
where $J$ denotes the matrix $\sum_{r\le 0}E_{rr}$, and for a
matrix $C=(c_{ij})\in\gltwo$, $\rm Str(C)$ stands for the
supertrace of the matrix $C$, which is given by
$\sum_{r\in\hf\Z}(-1)^{2r}c_{rr}$. We note that the expression
$\alpha(A,B)$ is well-defined for $A,B\in\gltwo$.  The restriction
of $\alpha$ to the subalgebras $\ospd$, $\pd_{\pm\pm}$ gives rise
to subalgebras of $\hgltwo$, which we denote accordingly by
$\widehat{\ospd}$ and $\widehat{\pd}_{\pm\pm}$, respectively.

Th linear map $\sigma$ sending $D$ to $D+1$ induces an associative
algebra automorphism of $\sd$.  Furthermore $(1-\sigma)$ maps
$\sd$ onto $\sd$, thus the map $(1-\sigma):\sd/{{\rm
ker}(1-\sigma)}\rightarrow\sd$ is a linear isomorphism.  Now let
${\rm Str}_s:\sd\rightarrow\C$, $s\in\C$, be the linear map
defined by $${\rm Str}_s(t^r\begin{pmatrix} f^0(D)&f^+(D)\\
f^-(D)&f^1(D))
\end{pmatrix}=\delta_{r,0}(f^0(s)-f^0(0)-f^1(s)+f^1(0)).$$
We note that ${\rm Str}_s$ vanishes on ${\rm ker}(1-\sigma)$ and
hence the linear map ${\rm
Str}_s\circ(1-\sigma)^{-1}:\sd\rightarrow\C$ is well-defined.  The
following lemma is a consequence of a straightforward computation
(cf. \cite{KR1}).

\begin{lem}\label{cocyclecomp}
For $t^rF(D),t^sG(D)\in\sd$, $r,s\in\Z$, we have
\begin{eqnarray*}
&&\Psi(t^rF(D),t^sG(D))+{\rm
Str}_s\circ(1-\sigma)^{-1}([t^rF(D),t^sG(D)]) \\ &=&
\alpha(\phi_s(t^rF(D)),\phi_s(t^sG(D))).
\end{eqnarray*}
\end{lem}

It follows from \lemref{cocyclecomp} that the map
$\widehat{\phi}_s:\hsd\rightarrow\hgltwo$ defined by
\begin{equation}
\widehat{\phi}_s(t^rF(D)) ={\phi}_s(t^rF(D)) +{\rm
Str}_s\circ(1-\sigma)^{-1}(t^rF(D))
\end{equation}
is a homomorphism of Lie superalgebras.

Let $x$ be a formal variable and define the generating functions
$e^{xD\theta\parth}$ and $e^{xD\parth\theta}$ in the usual way by
$\sum_{n\ge 0}\frac{(xD\theta\parth)^n}{n!}$ and $\sum_{n\ge
0}\frac{(xD\parth\theta)^n}{n!}$, respectively.  Since
$(1-\sigma)^{-1}(e^{xw})=-\frac{e^{xw}-1}{e^x-1}$ we have the
following description of the map $\widehat{\phi}_s$. A somewhat
different proof can be given parallel to the proof of Lemma 5.1,
pp.87, \cite{KWY}.

\begin{prop} \label{propext}
The $\C$-linear map $\widehat{\phi}_s:\hsd\rightarrow\hgltwo$
defined by
\begin{align*}
 \widehat{\phi}_s|_{\hsd_{j}} =& \phi_s|_{\hsd_{j}},\quad j\not=0,\\
 \widehat{\phi}_s(e^{xD\parth\theta})
 =& \phi_s(e^{xD\parth\theta})-\frac{e^{sx}-1}{e^x-1}\cdot {\bf 1},\\
 \widehat{\phi}_s(e^{xD\theta\parth})
 =& \phi_s(e^{xD\theta\parth})+\frac{e^{sx}-1}{e^x-1}\cdot {\bf 1},\\
 \widehat{\phi}_s(C)=& {\bf 1},
\end{align*}
is a homomorphism of Lie superalgebras.
\end{prop}

We note that it follows from \propref{sdzeroform} that
$\hsdzero_0$ is spanned by the basis elements
$\{D^n\parth\theta,(D+\hf)^n\theta\parth|n\ {\rm odd}\}$. We thus
obtain the following description for the embedding of $\hsdzero$
into $\hgltwo$ by restricting the map $\widehat{\phi}_s$. (Compare
Proposition~5.2, pp.87, \cite{KWY}).

\begin{prop}\label{propext2}
The $\C$-linear map $\widehat{\phi}_s:\hsdzero\rightarrow\hgltwo$
defined below is a homomorphism of Lie superalgebras:
\begin{align*}
 \widehat{\phi}_s|_{\hsdzero_{j}}
 =& \phi_s|_{\hsdzero_{j}},\quad j\not=0,\\
 \widehat{\phi}_s({\rm sinh}(xD\parth\theta))
 =& \phi_s({\rm sinh} (xD\parth\theta))
  -\frac{{\rm cosh}((s-\hf)x)
  -{\rm cosh}(\frac{x}{2})}{{\rm sinh}(\frac{x}{2})}\cdot {\bf 1},\\
  \widehat{\phi}_s({\rm sinh}(x(D+\hf)\theta\parth))
  =& \phi_s({\rm sinh}(x(D+\hf)\theta\parth))
  +\frac{{\rm cosh}(sx)-1}{{\rm sinh}(\frac{x}{2})}\cdot {\bf 1},\\
 \widehat{\phi}_s(C)=& {\bf 1}.
\end{align*}
In particular, $\widehat{\phi}_0$ is a homomorphism from
$\hsdzero$ to $\widehat{\ospd}$.
\end{prop}

It follows from \propref{zeroform1}, \propref{zeroform2},
\propref{zeroform3} and \propref{zeroform4}  that the Cartan
subalgebra ${}^{\pm\pm}\widehat{\sd}_0$ is spanned by the basis
elements
$\{(D+\hf)^n\parth\theta-(-1)^n(D+\hf)^n\theta\parth|n\in\Z_+\}$
which is the same as
$\{(D+\hf)^{2k+1},(D+\hf)^{2k}(\parth\theta-\theta\parth)|k\in\Z_+\}$.
We may consider the following generating functions:

\begin{equation*}
{\rm sinh}((D+1/2)x), \quad {\rm cosh}
((D+1/2)x)(\parth\theta-\theta\parth) .
\end{equation*}
Thus restricting $\hat{\psi}_s$ we obtain the following
description for the embedding of ${}^{\pm\pm}\widehat{\sd}$ into
$\hgltwo$.

\begin{prop}
The $\C$-linear map
$\widehat{\phi}_s:{}^{\pm\pm}\widehat{\sd}\rightarrow\hgltwo$
defined below is a homomorphism of Lie superalgebras:
\begin{align*}
 \widehat{\phi}_s|_{{}^{\pm\pm}\widehat{\sd}_{j}}
 &=\phi_s|_{{}^{\pm\pm}\widehat{\sd}_{j}},\quad
 j\not=0,\\
 \hat{\phi}_s({\rm sinh}((D+\hf)x))
 &=\phi_s({\rm sinh}((D+\hf)x)),\\
 \widehat{\phi}_s({\rm cosh}((D+\hf)x)(\parth\theta-\theta\parth))
 &= \phi_s({\rm cosh}((D+\hf)x)(\parth\theta-\theta\parth))\\
 &\qquad
 -2\frac{{\rm cosh}((s+\hf)x)-{\rm cosh}(\frac{x}{2})}{e^x-1}\cdot{\bf 1},\\
 \widehat{\phi}_s(C) &={\bf 1}.
\end{align*}
In particular, $\widehat{\phi}_0$ is a homomorphism from
${}^{\pm\pm}\widehat{\sd}$ to $\hpd_{\pm\pm}$.
\end{prop}

We conclude this section with the following proposition whose
proof is analogous to the one of Proposition 4.3 in \cite{KR1}.

\begin{prop}\label{holoext}
Let $\G$ be either $\hsd$ or ${}^i\hsd$, where $i=0,\pm\pm$.  Let
$\G^{\co}$ denote the corresponding holomorphic completion of $\G$
and $V$ a quasifinite $\G$-module.  The the action of $\G$ extends
naturally to an action of at least $\G^{\co}_{k}$, for all
$k\not=0$.
\end{prop}

\section{Finite-dimensional Howe duality} \label{sect_howedual}

Howe duality \cite{H1, H2} was generalized in \cite{CW1, CW2} in a
systematic way (also cf. \cite{OP, Naz, Se} for other independent
approaches\footnote{We thank M.~Nazarov for bringing the reference
\cite{Naz, OP} to our attention.}) to finite-dimensional Lie
superalgebras. However in order to apply those results to our
infinite dimensional setting, we need to use Borel subalgebras
different from the ones used in \cite{CW1}. In this section we
will obtain explicit formulas of the joint highest weight vectors
for a dual pair of Lie superalgebras acting on a supersymmetric
module with respect to the new Borel subalgebras, and thus obtain
an explicit isotypic decomposition of the supersymmetric module
with respect to the action of dual pairs. A variation of these
formulas will be applied to the study of Howe duality involving
infinite-dimensional Lie superalgebras in the next section.
%
%
%
%
%
%

Let $\gl(m|n)$ denote the general linear Lie superalgebra of
linear transformations on the $(m|n)$-dimensional complex
superspace  $\C^{m|n}$.  Let the subspace $\sum_{i=1}^{m+n}\C
E_{ii}$ of diagonal matrices be our choice of a Cartan subalgebra
$\h$ and $B$ any Borel subalgebra containing $\h$. A
finite-dimensional irreducible representation $V_{m|n}^\la$ is
determined by a unique (up to a scalar) non-zero highest weight
vector $v\in V_{m|n}^\la$, where $\la\in\h^*$, which in turn is
uniquely determined by the properties
\begin{eqnarray*}
h v&=&\la(h)v,\quad  h\in\h,\\ b v&=&0,\quad b\in[B,B].
\end{eqnarray*}

Consider the natural action of the Lie superalgebra $\gl(n|n)$ on
the $\C^{n|n}$ and the natural action of the Lie algebra $\gl(l)$
on the space $\C^l$.  We obtain an action of
$\gl(n|n)\times\gl(l)$ on the space $\C^{n|n}\otimes\C^l$.

We will {\bf assume for the rest of this section that $n\ge l$},
which is all we need for the application in the next section.

We consider the induced action of $\gl(n|n)\times\gl(l)$ on the
$k$-th skew-symmetric tensor $\Lambda^k(\C^{n|n}\otimes\C^l)$.  As
a $\gl(n|n)\times\gl(l)$-module $\La^k(\C^{n|n}\otimes\C^l)$ is
completely reducible and we have

\begin{equation}\label{glduality}
\Lambda^k(\C^{n|n}\otimes\C^l)\cong\sum_{\la}V^{\la'}_{n|n}\otimes
V^\la_l,
\end{equation}
where $\la=(\la_1,\la_2,\ldots )$ is summed over all partition of
size $|\la|=k$ and length $l(\la)\le l$.  As usual $\la'$ denotes
the transpose of $\la$. Here the partition $\la$ is viewed as
highest weights of $\gl(n|n)$ and $\gl(l)$ with respect to their
standard Borel subalgebras of upper triangular matrices as
follows. As a $\gl(l)$-highest weight we have $\la(E_{ii})=\la_i$,
for $i=1,\ldots,l$, while as a $\gl(n|n)$-highest weight we have
$\la(E_{ii})=\la'_i$, for $i=1,\ldots,n$, and $\la(E_{j+n,j+n})=
\langle \la_j-n\rangle$, for $j=1,\ldots,n$ (cf. \cite{CW1, OP,
Se}). The symbol $\langle m \rangle$, $m\in\Z$, stands for $m$, if
$m\in\Z_+$, and stands for zero otherwise.

For our application in next section, we will use the technique
developed in \cite{CW1} to derive explicit formulas for the joint
highest weight vectors with respect to certain non-standard Borel
subalgebras (compare \cite{OP} where less explicit formulas are
also presented).

We note that while all Borel subalgebras of $\gl(l)$ are conjugate
to each other, this is no longer true for Lie superalgebras
\cite{K}. Therefore it follows that a different choice of Borel
subalgebra for $\gl(n|n)$ above may result in different
description of $\gl(n|n)$-highest weights in the decomposition
\eqnref{glduality}.

Let $\{e^1,\ldots,e^n,e^{n+1},\ldots,e^{n+n}\}$ be the standard
basis of $\C^{n|n}$. In particular $\{e^i|1\le i\le n\}$ are even
basis vectors, while $\{e^i|n+1\le i\le 2n\}$ are odd basis
vectors. Of course any ordering of this basis gives rise, via
taking the upper triangular matrices with respect to this
ordering, to a Borel subalgebra of $\gl(n|n)$ containing the
Cartan subalgebra of diagonal matrices. We consider two orderings
that are relevant for our discussion. The first one is the
ordering $\{e^1,e^{n+1},e^2,e^{n+2},\ldots,e^{n},e^{n+n}\}$ and we
denote by $B_1$ the Borel subalgebra corresponding to this
ordering.  The second ordering is
$\{e^{n+1},e^1,e^{n+2},e^2,\ldots,e^{n+n},e^n\}$ and we denote the
corresponding Borel subalgebra by $B_2$. Our present goal is to
find the $\gl(n|n)\times\gl(l)$-joint highest weight vectors
inside $\La^k(\C^{n|n}\otimes\C^l)$ with respect to the Borel
subalgebra $B_i\times B'$, $i=1,2$, where $B'$ is the Borel
subalgebra of $\gl(l)$ consisting of upper triangular $l\times l$
matrices.

Denote by $e_j$, $j=1,\ldots,l$ the standard basis for $\C^l$.
Then $\{e^{i}\otimes e_j|1\le i\le 2n,1\le j\le l\}$ is a basis
for $\C^{n|n}\otimes \C^l$.  We set for $1\le i\le n$ and $1\le
j\le l$:
\begin{eqnarray*}
\xi^i_j&=&e^{i}\otimes e_j,\\ x^i_j&=&e^{n+i}\otimes e_j.
\end{eqnarray*}
Then the skew-symmetric algebra
$\La^*(\C^{n|n}\otimes\C^l)=\sum_{k=0}^\infty\La^k(\C^{n|n}\otimes\C^l)$
may be identified with the superalgebra
$\C[x^i_j]\otimes\La(\xi^i_j)$, the tensor algebra of the
polynomial algebra in the variables $x^i_j$ and the Grassmann
superalgebra in the variables $\xi^i_j$, for $1\le i\le n$ and
$1\le j\le l$.

Under this identification the action of $\gl(n|n)$  on
$\C[x^i_j]\otimes\La(\xi^i_j)$ may be realized as first order
linear differential operators
\begin{eqnarray*}
\sum_{j=1}^l \xi^i_j\frac{\pa}{\pa\xi^k_j},\quad\sum_{j=1}^l
x^i_j\frac{\pa}{\pa\xi^k_j},\quad\sum_{j=1}^l x^i_j\frac{\pa}{\pa
x^k_j},\quad\sum_{j=1}^l \xi^i_j\frac{\pa}{\pa x^k_j},\quad 1\le
i,k\le n,
\end{eqnarray*}
while that of $\gl(l)$ may be realized as
\begin{eqnarray*}
\sum_{i=1}^n\xi^i_j\frac{\pa}{\pa\xi^i_k}+ \sum_{i=1}^n
x^i_j\frac{\pa}{\pa x^i_k},\quad 1\le j,k\le l.
\end{eqnarray*}

The Cartan subalgebra of $\gl(n|n)$ is spanned by the basis
vectors $\sum_{j=1}^l\xi_j^i\frac{\pa}{\pa\xi^i_j}$ and
$\sum_{j=1}^l x^i_j\frac{\pa}{\pa x^i_j}$, for $1\le i\le n$,
while the nilpotent radical of the Borel subalgebra $B_1$ is
generated by the odd simple root vectors

\begin{eqnarray}
&\sum_{j=1}^l \xi^i_j\frac{\pa}{\pa x^i_j}, \quad 1\le i\le
n;\label{glnnrad1}\\
&\sum_{j=1}^l x^i_j\frac{\pa}{\pa\xi^{i+1}_j},\quad 1\le i\le
n-1.\label{glnnrad2}
\end{eqnarray}
The nilpotent radical of the Borel subalgebra $B_2$ is generated
by the odd simple root vectors
\begin{eqnarray}
&\sum_{j=1}^l x^i_j\frac{\pa}{\pa \xi^i_j}, \quad 1\le i\le
n;\label{glnnrad3}\\
&\sum_{j=1}^l \xi^i_j\frac{\pa}{\pa x^{i+1}_j},\quad i\le i\le
n-1.\label{glnnrad4}
\end{eqnarray}

The Cartan subalgebra of $\gl(l)$ on the other hand is spanned by
$\sum_{i=1}^n\xi^i_j\frac{\pa}{\pa\xi^i_j}+ \sum_{i=1}^n
x^i_j\frac{\pa}{\pa x^i_j}$ , for $1\le j\le l$, while the
nilpotent radical of the Borel subalgebra is generated by the
simple root vectors
\begin{align}\label{glrad}
\sum_{i=1}^n\xi^i_j\frac{\pa}{\pa\xi^i_{j+1}}+ \sum_{i=1}^n
x^i_j\frac{\pa}{\pa x^i_{j+1}},\quad 1\le j\le l-1.
\end{align}

Consider the following $l\times l$ matrices

\begin{eqnarray*}
X^1=&
\begin{pmatrix}
\xi_1^1&\xi_2^1&\cdots &\xi^1_l\\ \xi_1^1&\xi_2^1&\cdots
&\xi^1_l\\ \vdots&\vdots&\cdots &\vdots\\ \xi_1^1&\xi_2^1&\cdots
&\xi^1_l\\
\end{pmatrix},\allowdisplaybreaks\\
X^2=&
\begin{pmatrix}
x_1^1&x_2^1&\cdots &x^1_l\\ \xi_1^2&\xi_2^2&\cdots &\xi_l^2\\
\vdots&\vdots&\cdots &\vdots\\ \xi_1^2&\xi_2^2&\cdots &\xi^2_l\\
\end{pmatrix},\allowdisplaybreaks\\
X^3=&
\begin{pmatrix}
x_1^1&x_2^1&\cdots &x_l^1\\ x_1^2&x_2^2&\cdots &x_l^2\\
\xi_1^3&\xi_2^3&\cdots &\xi^3_l\\ \vdots&\vdots&\cdots &\vdots\\
\xi_1^3&\xi_2^3&\cdots &\xi^3_l\\
\end{pmatrix},\allowdisplaybreaks\\
&\vdots\\ X^k \equiv X^{l-1}=&
\begin{pmatrix}
x_1^1&x_2^1&\cdots &x_l^1\\ x_1^2&x_2^2&\cdots &x_l^2\\
x_1^3&x_2^3&\cdots &x^3_l\\ \vdots&\vdots&\cdots &\vdots\\
x_1^l&x_2^l&\cdots &x_l^l\\
\end{pmatrix},\quad k \ge l.\allowdisplaybreaks\\
\end{eqnarray*}

For $0\le r\le l$ we let $X^i_r$ denote the first $r\times r$
minor of the matrix $X^i$.

\begin{rem}
By the determinant of an $r\times r$ matrix $A=(a_i^j)$, denoted
by ${\rm det}A$, we will always mean the expression
$\sum_{\sigma\in S_r}{\rm
sgn}(\sigma)a^1_{\sigma_1}a^2_{\sigma_2}\cdots a^r_{\sigma_r}$. It
is necessary to specify the order of product as some of the
$a^j_i$'s might be odd variables.
\end{rem}

The following lemma, which is the Corollary~4.1 in \cite{CW1},
plays an important role.

\begin{lem}\label{identity}
Let $x_i^j$ be even variables for $i=1,\ldots,r$ and
$j=1,\ldots,q$ with $r\ge t>q$.  Let $\xi^q_i$ and $\xi^s_i$ be
odd variables for $i=1,\ldots, l$.  Then $${\rm det}
\begin{pmatrix}
x^1_1&x^1_2&\cdots&x^1_r\\ x^2_1&x_2^2&\cdots&x^2_r\\
\vdots&\vdots&\cdots&\vdots\\ x^q_1&x^q_2&\cdots&x^q_r\\
\xi^q_1&\xi^q_2&\cdots&\xi^q_r\\ \xi^q_1&\xi^q_2&\cdots&\xi^q_r\\
\vdots&\vdots&\cdots&\vdots\\ \vdots&\vdots&\cdots&\vdots\\
\xi^q_1&\xi^q_2&\cdots&\xi^q_r\\
\end{pmatrix}
\cdot{\rm det}
\begin{pmatrix}
x^1_1&x^1_2&\cdots&x^1_t\\ x^2_1&x_2^2&\cdots&x^2_t\\
\vdots&\vdots&\cdots&\vdots\\ x^q_1&x^q_2&\cdots&x^q_t\\
\xi^q_1&\xi^q_2&\cdots&\xi^q_t\\ \xi^s_1&\xi^s_2&\cdots&\xi^s_t\\
\vdots&\vdots&\cdots&\vdots\\ \xi^s_1&\xi^s_2&\cdots&\xi^s_t\\
\end{pmatrix}=0.
$$
\end{lem}

\begin{cor}\label{keycor}
Let $x_i^j$ be even variables for $i=1,\ldots,r$ and
$j=1,\ldots,q+m$ with $r\ge t>q$. Let $\xi^q_i$ and $\xi^s_i$ be
odd variables for $i=1,\ldots, l$. Then $${\rm det}
\begin{pmatrix}
x^1_1&x^1_2&\cdots&x^1_r\\ x^2_1&x_2^2&\cdots&x^2_r\\
\vdots&\vdots&\cdots&\vdots\\ x^q_1&x^q_2&\cdots&x^q_r\\
\xi^q_1&\xi^q_2&\cdots&\xi^q_r\\ \xi^q_1&\xi^q_2&\cdots&\xi^q_r\\
\vdots&\vdots&\cdots&\vdots\\ \vdots&\vdots&\cdots&\vdots\\
\xi^q_1&\xi^q_2&\cdots&\xi^q_r\\
\end{pmatrix}
\cdot{\rm det}
\begin{pmatrix}
x^1_1&x^1_2&\cdots&x^1_t\\ \vdots&\vdots&\cdots&\vdots\\
x^q_1&x^q_2&\cdots&x^q_t\\ \xi^q_1&\xi^q_2&\cdots&\xi^q_t\\
x^{q+2}_1&x^{q+2}_2&\cdots&x^{q+2}_t\\
\vdots&\vdots&\cdots&\vdots\\
x^{q+m}_1&x^{q+m}_2&\cdots&x^{q+m}_t\\
\xi^s_1&\xi^s_2&\cdots&\xi^s_t\\ \vdots&\vdots&\cdots&\vdots\\
\xi^s_1&\xi^s_2&\cdots&\xi^s_t\\
\end{pmatrix}=0.
$$
\end{cor}

\begin{proof}
The identity follows by applying successively the differential
operators of the form $\sum_{i=1}^t
x_i^{q+2}\frac{\partial}{\partial\xi_i^s}$, $\sum_{i=1}^t
x_i^{q+3}\frac{\partial}{\partial\xi_i^s}$,$\cdots$, $\sum_{i=1}^t
x_i^{q+m}\frac{\partial}{\partial\xi_i^s}$ to the identity in
\lemref{identity}.
\end{proof}

\begin{thm}\label{hwv1}
Let $\la$ be a partition of length $l(\la)\le l$.  Then the
expression
\begin{equation*}
\prod_{i=1}^{\la_1}{\rm det}X^i_{\la_i'}
\end{equation*}
is annihilated by \eqnref{glnnrad1}, \eqnref{glnnrad2} and
\eqnref{glrad}, and hence is a $\gl(n|n)\times\gl(l)$ highest
weight vector with respect to the Borel subalgebra $B_1\times B'$.
\end{thm}

\begin{proof}
From the definition of the determinant it is easy to see that the
simple root vectors in \eqnref{glrad} kill each ${\rm
det}X^i_{\la_i'}$. Hence they kill $\prod_{i=1}^{\la_1}{\rm
det}X^i_{\la_i'}$.

Also it is clear that $\prod_{i=1}^{\la_1}{\rm det}X^i_{\la_i'}$
is annihilated by \eqnref{glnnrad2}, since such an expression acts
on ${\rm det}X^{l+1}_r$ by replacing the row
$(\xi^{l+1}_1,\xi^{l+1}_2,\ldots,\xi^{l+1}_r)$ by
$(x^{l}_1,x^{l}_2,\ldots,x^{l}_r)$, where $l \geq 1$.  As the
latter row already appears in $X^{l+1}_r$, the resulting
determinant is zero.

Hence it remains to show that \eqnref{glnnrad1} also annihilates
$\prod_{i=1}^{\la_1}{\rm det}X^i_{\la_i'}$.  But this now can be
easily seen using \corref{keycor}.
\end{proof}

Let $\la$ be a partition of size $k$ and length $l(\la)\le l$.
Define $\phi_1(\la)$ to be the $n|n$-tuple of integers
$$\phi_1(\la) :=(\la_1',\langle\la_2'-1\rangle, \ldots,
\langle\la_n'-n+1\rangle; \langle\la_1-1\rangle,
\langle\la_2-2\rangle,\ldots, \langle\la_n-n\rangle).$$ The
following corollary follows from an easy computation of the weight
of vector in \thmref{hwv1}.

\begin{cor} As a $\gl(n|n)\times\gl(l)$-module we have
\begin{equation*}
\Lambda^k(\C^{n|n}\otimes\C^l)\cong\sum_{\la}V^{\phi_1(\la)}_{n|n}\otimes
V^\la_l,
\end{equation*}
where $\phi_1(\la)$ is the highest weight of $\gl(n|n)$ with
respect to the Borel subalgebra $B_1$ and $\la$ above is summed
over all partitions with $l(\la)\le l$ and $|\la|=k$.
\end{cor}



\begin{thm}\label{hwv2}
Let $\la$ be a partition of length $l(\la)\le l$.  Then the
expression
\begin{equation*}
\prod_{j=1}^{\la_1}{\rm det}X^{j+1}_{\la_j'}
\end{equation*}
is annihilated by \eqnref{glnnrad3}, \eqnref{glnnrad4} and
\eqnref{glrad}, and hence is a $\gl(n|n)\times\gl(l)$ highest
weight vector with respect to the Borel subalgebra $B_2\times B'$.
\end{thm}

\begin{proof}
The proof is similar to that of \thmref{hwv1}. Again it is easy to
show that the expression is annihilated by operators of the form
\eqnref{glrad} and \eqnref{glnnrad3}.  The fact that it is
annihilated by operators of the form \eqnref{glnnrad4} is again a
consequence of \corref{keycor}.
\end{proof}

Similarly for a partition $\la$ with $|\la|=k$ and $l(\la)\le l$
we define an $n|n$-tuple of integers
$$\phi_2(\la):=(\langle\la_1'-1\rangle, \langle
\la_2'-2\rangle,\ldots, \langle \la_n'-n\rangle; \la_1, \langle
\la_2-1\rangle,\ldots, \langle\la_n-n+1\rangle).$$

\begin{cor} As a $\gl(n|n)\times\gl(l)$-module we have
\begin{equation*}
\Lambda^k(\C^{n|n}\otimes\C^l)\cong\sum_{\la}V^{\phi_2(\la)}_{n|n}\otimes
V^\la_l,
\end{equation*}
where $\phi_2(\la)$ is the highest weight of $\gl(n|n)$ with
respect to the Borel subalgebra $B_2$ and $\la$ above is summed
over all partitions with $l(\la)\le l$ and $|\la|=k$.
\end{cor}

\begin{rem}
The arguments and results in this section remain valid for the
dual pair ${\rm gl}(m|n)\times {\rm gl}(l)$, as long as $m,n\ge
l$.  In this more general setting we need to replace $B_1$ and
$B_2$ with appropriate Borel subalgebras.  Let's denote the
standard basis of $\C^{m|n}$ by
$\{e^1,\cdots,e^m;e^{m+1},\cdots,e^{m+n}\}$. We may replace $B_1$
with any Borel subalgebra corresponding to any ordering of the
basis according to which
$e^{m+1},e^1,e^{m+2},e^2,\cdots,e^{m+l},e^l$ appear as the first
$2l$ members.  Similarly we may replace $B_2$ by any Borel
subalgebra whose corresponding order of basis has
$e^1,e^{m+1},e^2,e^{m+2},\cdots,e^l,e^{m+l}$ as its first $2l$
members.
\end{rem}

\begin{rem}
One can modify some results in the next section to obtain a finite
dimensional $({\rm Pin}(2l), {\rm osp}(2n+ 1,2m))$ duality ($n, m
\ge 1$) which seems to be new as well.
\end{rem}
\section{Free field realizations and duality}
\label{sect_freefield}

Take a pair of fermionic fields ($bc$ fields) $$
 \psi^{+}(z) = \sum_{n \in {\Z} } \psi^{+}_n z^{ -n -1},
  \quad \psi^{-}(z) = \sum_{n \in {\Z} }
                           \psi^{-}_n z^{-n},
$$ with the following anti-commutation relations
\begin{eqnarray*}
 [\psi^{+}_m, \psi^{-}_n ]  =  \delta_{m+n, 0}, \quad
 {[} {\psi}^{\pm}_m, {\psi}^{\pm}_n {]}  =  0,\quad m, n \in \Z.
\end{eqnarray*}

Also we take a pair of bosonic ghost fields ($\beta \gamma$
fields) $$ \gamma^{\pm} (z) = \sum_{ r \in \hz }
                         \gamma^{\pm}_r z^{-r - \hf} ,
$$ with the following commutation relations
\begin{eqnarray*}
 [ \gamma^{+}_r, \gamma^{-}_s ] =\delta_{r+s,0}, \quad
 [ \gamma^{\pm}_r, \gamma^{\pm}_s ] =0,\quad r, s \in \hf + \Z.
\end{eqnarray*}

Denote by $\F$ the Fock space of $\psi^{\pm }(z)$ and $\gamma^{\pm
}(z)$ generated by a vacuum vector $\vac$ which is annihilated by
$\psi^{+}_m, \psi^{-}_{m +1} (m \in \Z_+) $, $ \gamma^{\pm}_r$
$(r\in \hf + \Z_+)$.

More generally we take $l$ (independent) copies of $bc\beta
\gamma$ fields $\psi^{\pm,k} (z)$,  $\gamma^{\pm,k} (z) $ $(k = 1,
\dots, l)$ and consider the corresponding Fock space $\F^{\otimes
l}$.
\subsection{The case of $\hgltwo$ and $\hsd$}

In this subsection, we will realize $\hgltwo$ and $\hsd$ in the
Fock space $\F^{\otimes l}$, also cf. \cite{AFMO1}. We then
establish a $(GL(l), \gltwo)$ duality in $\F^{\otimes l}$. A
$(GL(l), \hsd)$ duality follows from this via the embeddings
$\phi_s$.

Introduce the following generating functions for $\hgltwo$:
\begin{eqnarray*}
E_0 (z,w) &=&\sum_{i,j \in \Z} E_{ij} z^{i -1} w^{ -j},
\\
E_1 (z,w) &=&\sum_{r,s \in\hf+\Z} E_{rs} z^{r-\hf} w^{-s-\hf},
\\
E_+ (z,w) &=&\sum_{i\in \Z,s\in\hf+\Z} E_{is} z^{i-1} w^{ -s-\hf},
\\
E_- (z,w) &=&\sum_{r\in\hf+\Z,j\in\Z} E_{rj} z^{r-\hf}w^{-j}.
\end{eqnarray*}

We have the following free field realization for $\hsd$.

\begin{prop}
  Let
 \begin{eqnarray*}
   E_0 (z, w) & =& \sum_{k=1}^l :\psi^{+,k} (z) \psi^{-,k} (w): , \\
   E_1 (z, w) & =& - \sum_{k=1}^l :\gamma^{+,k} (z) \gamma^{-,k} (w): , \\
   E_+ (z, w) & =& \sum_{k=1}^l :\psi^{+,k} (z) \gamma^{-,k} (w):  , \\
   E_- (z, w) & =&-  \sum_{k=1}^l :\gamma^{+,k} (z) \psi^{-,k} (w): .
 \end{eqnarray*}
 This defines a representation of $\hgltwo$ on $\F^{\otimes l}$ of central charge $l$.
\end{prop}

Componentwise, we can write down the above representation of
$\hgltwo$ on $\F^{\otimes l}$ as follows:

\begin{eqnarray*}
 E_{ij} = \sum_{k=1}^l : \psi^{+,k}_{-i} \psi^{-,k}_j:, &&
 E_{rs} = -\sum_{k=1}^l : \gamma^{+,k}_{-r} \gamma^{-,k}_s:, \\
 E_{is} = \sum_{k=1}^l : \psi^{+,k}_{-i} \gamma^{-,k}_s:, &&
 E_{rj} = -\sum_{k=1}^l : \gamma^{+,k}_{-r} \psi^{-,k}_j:,
\end{eqnarray*}
where $i,j \in \Z$ and $r,s \in \hf+ \Z$. The normal ordering
$:\;:$ is defined to move the annihilation operators to the right.

It is well known that the fields $:\psi^{+,i}(z)\psi^{j,-}(z):$
define a representation of the affine algebra $\widehat{\rm
gl}(l)$ on the Fock space $\F^{\otimes l}$ of central charge $l$.
On the other hand the components of fields
$-:\gamma^{+,i}(z)\gamma^{-,j}(z):$ define a representation of
$\widehat{\rm gl}(l)$ on $\F^{\otimes l}$ of central charge $-l$.
Hence
\begin{equation*}
:\psi^{+,i}(z)\psi^{j,-}(z):-:\gamma^{+,i}(z)\gamma^{-,j}(z):,\quad
1\le i,j\le l,
\end{equation*}
give a representation of $\widehat{\rm gl}(l)$ of central charge
zero on $\F^{\otimes l}$. The horizontal subalgebra is ${\rm
gl}(l)$ and its action on $\F^{\otimes l}$ lifts to a rational
action of the group ${\rm GL}(l)$.  In particular as a ${\rm
GL}(l)$-module $\F^{\otimes l}$ is completely reducible. The
following proposition is a variant in our infinite-dimensional
setting of the description of ${\rm GL}(l)$-invariants in ${\rm
End}(\F^{\otimes l})$ \cite{H1}.

\begin{prop} \label{propfree}
${\rm GL}(l)$ and $\hgltwo$ are mutual centralizers in ${\rm End
}(\F^{\otimes l})$.  In particular $\F^{\otimes l}$ is a
multiplicity-free direct sum of irreducible $\hgltwo\times{\rm
GL}(l)$-modules.
\end{prop}

Our next task is to determine the decomposition of $\F^{\otimes
l}$ with respect to the joint action of $\hgltwo\times{\rm
GL}(l)$.

We will do this by explicitly displaying all joint
$\hgltwo\times{\rm GL}(l)$ highest weight vectors with respect to
the (joint) Borel subalgebra $\hat{B}_{\infty|\infty}\times B'$,
where $\hat{B}_{\infty|\infty}=\oplus_{r\ge 0}(\hgltwo)_r$.

For this purpose we need some notation.  Define for $j\in\Z_+$ the
matrices $X^{-j}$ as follows:
\begin{eqnarray*}
 X^0=&
\begin{pmatrix}
\psi_0^{-,l}&\psi_0^{-,l-1}&\cdots &\psi_0^{-,1}\\
\psi_0^{-,l}&\psi_0^{-,l-1}&\cdots &\psi_0^{-,1}\\
\vdots&\vdots&\cdots &\vdots\\
\psi_0^{-,l}&\psi_0^{-,l-1}&\cdots
&\psi_0^{-,1}\\
\end{pmatrix},\allowdisplaybreaks\\
 X^{-1}=&
\begin{pmatrix}
\gamma_{-\hf}^{-,l}&\gamma_{-\hf}^{-,l-1}&\cdots
&\gamma_{-\hf}^{-,1}\\ \psi_{-1}^{-,l}&\psi_{-1}^{-,l-1}&\cdots
&\psi_{-1}^{-,1}\\ \vdots&\vdots&\cdots &\vdots\\
\psi_{-1}^{-,l}&\psi_{-1}^{-,l-1}&\cdots &\psi_{-1}^{-,1}\\
\end{pmatrix},\allowdisplaybreaks\\
 X^{-2}=&
\begin{pmatrix}
\gamma_{-\hf}^{-,l}&\gamma_{-\hf}^{-,l-1}&\cdots
&\gamma_{-\hf}^{-,1}\\
\gamma_{-\frac{3}{2}}^{-,l}&\gamma_{-\frac{3}{2}}^{-,l-1}&\cdots
&\gamma_{-\frac{3}{2}}^{-,1}\\
\psi_{-2}^{-,l}&\psi_{-2}^{-,l-1}&\cdots &\psi_{-2}^{-,1}\\
\vdots&\vdots&\cdots &\vdots\\
\psi_{-2}^{-,l}&\psi_{-2}^{-,l-1}&\cdots &\psi_{-2}^{-,1}\\
\end{pmatrix},\allowdisplaybreaks\\
&\vdots\\&\vdots\\
 X^{-k} \equiv {}^- X^{-l}=&
\begin{pmatrix}
\gamma_{-\hf}^{-,l}&\gamma_{-\hf}^{-,l-1}&\cdots
&\gamma_{-\hf}^{-,1}\\
\gamma_{-\frac{3}{2}}^{-,l}&\gamma_{-\frac{3}{2}}^{-,l-1}&\cdots
&\gamma_{-\frac{3}{2}}^{-,1}\\ \vdots&\vdots&\cdots &\vdots\\
\gamma_{-l+\hf}^{-,l}&\gamma_{-l+\hf}^{-,l-1}&\cdots
&\gamma_{-l+\hf}^{-,1}\\
\gamma_{-l-\hf}^{-,l}&\gamma_{-l-\hf}^{-,l-1}&\cdots
&\gamma_{-l-\hf}^{-,1}\\
\end{pmatrix}, \quad k \ge l.\allowdisplaybreaks\\
\end{eqnarray*}
The matrices $ X^j$, for $j\in\N$, are defined in a similar
fashion. Namely, $ X^j$ is obtained from $ X^{-j}$ by replacing
$\psi^{-,k}_i$ by $\psi^{+,l-k+1}_i$ and $\gamma^{-,k}_r$ by
$\gamma^{+,l-k+1}_r$.  For example,
\begin{eqnarray*}
  X^1=&
\begin{pmatrix}
\gamma_{-\hf}^{+,1}&\gamma_{-\hf}^{+,2}&\cdots
&\gamma_{-\hf}^{+,l}\\ \psi_{-1}^{+,1}&\psi_{-1}^{+,2}&\cdots
&\psi_{-1}^{+,l}\\ \vdots&\vdots&\cdots &\vdots\\
\psi_{-1}^{+,1}&\psi_{-1}^{+,2}&\cdots &\psi_{-1}^{+,l}\\
\end{pmatrix}
\end{eqnarray*}

For $0\le r\le l$, we let $X^i_{r}$($i \ge 0$) denote the first
$r\times r$ minor of the matrix $ X^i$ and let $X^i_{-r}$($i <0$)
denote the first $r\times r$ minor of the matrix $ X^i$.

Consider a generalized Young diagram
$\la=(\la_1,\la_2,\cdots,\la_p)$ of length $l$ with
\begin{equation*}
\la_1\ge\la_2\ge\cdots\ge\la_i>\la_{i+1}
=0=\cdots=\la_{j-1}>\la_j\ge\cdots\ge\la_l.
\end{equation*}
It is well-known that the irreducible rational representations of
${\rm GL}(l)$ are parameterized by generalized Young diagrams,
hence they may be interpreted as highest weights of irreducible
representations of ${\rm GL}(l)$.  As usual we denote by
$\la_{j}'$ the length of the $j$-th column of $\la$. We use the
convention that the first column of $\la$ is the first column of
the Young diagram $\la_1\ge\la_2\ge\cdots\ge\la_i$.  The column to
the right is the second column of $\la$, while the column to the
left of it is the zeroth column and the column to the left of the
zeroth column is the $-1$-st column.  We also use the convention
that a non-positive column has negative length. As an example
consider $\la=(5,3,2,1,-1,-2)$ with $l(\la)=6$. We have
$\la_{-1}'=-1$, $\la_0'=-2$ $\la_1'=4$ etc.~(see
\eqnref{standard}).
\begin{equation}\label{standard}
{\beginpicture \setcoordinatesystem units <1.5pc,1.5pc> point at 0
2 \setplotarea x from -1.5 to 1.5, y from -2 to 4 \plot 0 0 0 4 1
4 1 0 0 0 / \plot 1 1 2 1 2 4 1 4 / \plot 2 2 3 2 3 4 2 4 / \plot
3 3 5 3 5 4 3 4 / \plot 0 1 1 1 / \plot 0 2 2 2 / \plot 0 3 3 3 /
\plot 4 3 4 4 / \plot 0 0 -1 0 / \plot 0 0 0 -2 / \plot -1 0 -1 -2
/ \plot 0 -2 -2 -2 / \plot 0 -1 -2 -1 / \plot -2 -2 -2 -1 /
\endpicture}
\end{equation}

Let $\La\in (\hgltwo)_0^*$.  Then $\La$ may be interpreted as a
highest weight for a highest weight irreducible representation of
$\hgltwo$.  We let $\la_a=\La(E_{aa})$, for $a\in\hf\Z$. Given a
generalized Young diagram $\la$ with $l(\la)\le l$, we define
$\Lambda(\la)$ of $(\hgltwo)_0^*$ to be the weight whose
components are given by:

\begin{align*}
 \la_i&= \langle\la_i'-i \rangle,\quad i\in\N,\\
 \la_j&= -\langle-\la'_j+j \rangle, \quad j\in -\Z_{+},\\
 \la_r&=\langle\la_{r+1/2}-(r-1/2) \rangle,\quad r\in\hf+\Z_+,\\
 \la_{s}&=-\langle -\la_{p+(s+1/2)}+({s-1/2}) \rangle,\quad s\in-\hf-\Z_+.
\end{align*}

\begin{thm}\label{gltwohwv}
 \begin{enumerate}
  \item
As a ${\rm GL}(l)\times\hgltwo$-module, $\F^{\otimes l}$ is
completely reducible and decomposes into isotypic components as
follows:
\begin{equation*}
\F^{\otimes l}\cong\bigoplus_{\la} V_l^\la\otimes
L(\hgltwo,\Lambda(\la)),
\end{equation*}
where the summation is summed over all generalized Young diagrams
$\la$ with $l(\la)\le l$.
\item
The  ${\rm GL}(l)\times\hgltwo$-highest weight vector with respect
to the Borel subalgebra $\hat{B}_{\infty|\infty}\times B'$ in the
$\la$-isotypic component of $\F^{\otimes l}$ is
\begin{equation}
{\rm det} X^{\la_{l}+1}_{{\la}_{\la_{l-1}}'}\cdots{\rm det}
X^{-1}_{{\la}_{-1}'} \cdot {\rm det} X^0_{\la_0'} \cdot {\rm det}
X^1_{\la_1'} \cdot {\rm det} X^2_{\la_2'}\cdots{\rm det}
X^{\la_1}_{{\la}_{\la_1}'}\vac.
\end{equation}
\end{enumerate}
\end{thm}

\begin{proof}
The second part is a variation in our infinite-dimensional setting
of \thmref{hwv1}. Since $\F^{\otimes l}$ is a rational
representation of ${\rm GL}(l)$, it follows from the second part
that all irreducible representations of ${\rm GL}(l)$ appears in
the decomposition of $\F^{\otimes l}$. The multiplicity-freeness
follows from \propref{propfree}.
\end{proof}

The irreducible representation $L(\hgltwo, \Lambda(\la))$ pulls
back via $\hat{\phi}_s$ to an irreducible representation of
$\hsd$.  For $\xi \in (\hsd)_0^*$, we denote as usual by
$L(\hsd,\xi)$ the irreducible highest weight representation of
$\hsd$.  Introduce the generating functions
\begin{align*}
\Delta_0(x)&=-\xi(e^{xD\parth\theta}),\\
\Delta_1(x)&=-\xi(e^{xD\theta\parth}).
\end{align*}

\begin{prop}\label{hsdprop1}
Let $L(\hgltwo,\Lambda)$ be the irreducible highest weight
representation of highest weight $\Lambda$ and central charge $c$
such that $\Lambda(E_{aa}) = \lambda_a$. Then via $\hat{\phi}_s$
the module $L(\hgltwo,\La)$ pulls back to the irreducible highest
weight module $L(\hsd,\la_{\hsd})$, where $\la_{\hsd}$ is
specified by the generating functions
\begin{align*}
\Delta_0(x)&=
\frac{\sum_{i\in\Z}(\la_{i,i}-\la_{i+1,i+1}+\delta_{i,0}c)
e^{(-i+s)x}+c}{e^x-1},\\
 \Delta_1(x)&=
\frac{\sum_{i\in\Z}(\la_{i-\hf,i-\hf}-\la_{i+\hf,i+\hf}
-\delta_{i,0}c)e^{(-i+s)x}+c}{e^x-1}.
\end{align*}
\end{prop}

\begin{proof}
The generating functions above can be directly calculated by using
the embedding $\widehat{\phi}_s$ given in (\ref{eqembed}) and
\propref{propext}. The irreducibility follows from
\propref{propexactseq} and \propref{holoext}.
\end{proof}

\begin{rem} \label{remdualhsd}
It follows that $\hsd$ and ${\rm GL}(l)$ form a dual pair in the
sense of Howe on $\F^{\otimes l}$.  The decomposition of
$\F^{\otimes l}$ with respect to the joint action $\hsd\times{\rm
GL}(l)$ follows from \thmref{gltwohwv} and \propref{hsdprop1}.
\end{rem}

Indeed we can write down the action of $\hsd$ on $\F^{\otimes l}$
in an explicit manner. Note that $J^{a,k}_n \equiv J^k_nM_a$
($a=0,1,\pm$) together with $C$ span $\hsd$. Define

$$ J^{a,n}(z) =\sum_{k \in\Z} J_k^{a,n} z^{-n-k-1}.$$

\begin{prop}  \label{propfreehsd}
The action of $\hsd$ on $\F^{\otimes l}$ of central charge $l$ is
given in terms of the following generating functions:
\begin{eqnarray*}
 J^{0,n}(z) &=& -\sum_{p =1}^l :\psi^{+p}(z) \partial^n
 \psi^{-p}(z): \\
 J^{1,n}(z) &=& \sum_{p =1}^l :\gamma^{+p}(z) \partial^n
 \gamma^{-p}(z): \\
 J^{+,n}(z) &=& -\sum_{p =1}^l :\psi^{+p}(z) \partial^n
 \gamma^{-p}(z): \\
 J^{-,n}(z) &=& \sum_{p =1}^l :\gamma^{+p}(z) \partial^n
 \psi^{-p}(z):
\end{eqnarray*}
Here and further $\partial^n(\cdot)$ denotes the $n$-th derivative
with respect to $z$.
\end{prop}

\begin{proof}
We will prove the last identity only and the proof of the others
is similar. Using the embedding $\phi_0$ we calculate
  \begin{eqnarray*}
 J^{-,n}(z) &=&\sum_{k \in\Z} (-t^k[D]_n M_-) z^{-n-k-1} \\
        & =&\sum_{k \in \Z} -[-j]_n E_{j -k, j} M_-  z^{-k -n -1}   \\
        & =&\sum_{k \in \Z}  [-j]_n
        \sum_{p =1}^l :\gamma^{+p}_{k -j} \psi^{-p}_j : z^{-k -n -1}   \\
       & =& \sum_{p =1}^l
             :\gamma^{+p} (z)\partial^n \psi^{-p} (z):.
\end{eqnarray*}
\end{proof}

Let ${\mathfrak P} = \{ J^{a,n}_k \mid n +k \geq 0, k \in \Z, n
\in \Z_+, a=0,1,\pm \}$ and let $\widehat{\mathfrak P} =
{\mathfrak P} \bigoplus \C C.$ One can check that
$\widehat{\mathfrak P}$ is a parabolic subalgebra of $\hsd$ (the
central extension when restricted to $\mathfrak P$ is trivial).
Geometrically, ${\mathfrak P}$ consists of those differential
operators in $\sd$ which extend to the interior of the circle.
Denote by $M_c(\hsd)$ ($c\in\C$) the vacuum module (which is a
generalized Verma module) $$
  M(\hsd,\widehat{\mathfrak P}, \Lambda) = {\mathcal U} (\hsd)
                 \bigotimes_{ {\mathcal U}(\widehat{\mathfrak P})}
                 \C_c
$$ where ${\mathcal U}(\cdot)$ denotes the universal enveloping
algebra, $\C_c$ is the one-dimensional representation of
$\widehat{\mathfrak P}$ by letting $C = c {\rm Id}$ and
$\widehat{\mathfrak P}\cdot \C_l =0.$ Denote by $V_c(\hsd)$ the
irreducible quotient of the $\hsd$-module $M_c(\hsd)$. Then
$M_c(\hsd)$ and $V_c(\hsd)$ carry vertex superalgebra structures.
For example, a proof of this can be given based on the above free
field realization given in \propref{propfreehsd} parallel to the
proof of Theorem~14.1, pp. 132, \cite{KWY}.

Now the duality in \thmref{gltwohwv} (cf.~\remref{remdualhsd}) can
be interpreted as a duality between $GL(l)$ and the vertex
superalgebra $V_l(\hsd)$ on the Fock space $\F^{\otimes l}$. The
irreducible $\hgltwo$-module appearing in $\F^{\otimes l}$ becomes
irreducible module over the vertex superalgebra $V_l(\hsd)$.
\subsection{The case of $\widehat{\ospd}$ and $\hsdzero$}

In this subsection we will construct actions of $\hospd$ and
$\hsdzero$ on $\F^{\otimes l}$, and establish a $(\hospd, {\rm
Pin}(2l))$ duality and then a $(\hsdzero, {\rm Pin}(2l))$ duality.
Let
\begin{eqnarray*}
   &&   E_0 (z, w) + z^{-1} wE_0 (-w, -z)\\
    & = & \sum_{i,j \in \Z}(E_{ij} -(-1)^{i +j} E_{-j, -i}) M_0 \;
             z^{i-1} w^{ -j}\allowdisplaybreaks \\
    & = & \sum_{k=1}^l \left(  :\psi^{+,k} (z) \psi^{-,k} (w):
                              -z^{-1}w :\psi^{+,k} (-w) \psi^{-,k} ( -z):
                        \right) , \allowdisplaybreaks\\
   &&E_1 (z, w) + E_1 (-w, -z)\\
    & = & \sum_{i,j \in \Z} (E_{ij} -(-1)^{i +j}E_{1-j,1-i})M_1\;
             z^{i -1} w^{ -j}  \allowdisplaybreaks\\
    & = & -\sum_{k=1}^l \left(:\gamma^{+,k} (z) \gamma^{-,k} (w):
                              + :\gamma^{+k} (-,w) \gamma^{-,k} ( -z):
                        \right), \allowdisplaybreaks\\
   && z E_+ (z, w) +  E_- (-w, -z)\\
    & = & \sum_{i,j \in \Z} (E_{ij} M_+ +(-1)^{i +j} E_{1-j, -i} M_-)
             z^{i} w^{ -j}   \allowdisplaybreaks\\
    & = & \sum_{k=1}^l \left( z :\psi^{+,k} (z) \gamma^{-,k} (w):
                           - :\gamma^{+l} (-,w) \psi^{-,k} ( -z):
                        \right).
\end{eqnarray*}

The last equation above is equivalent to
 $$E_- (z, w) -wE_+ (-w, -z) = -\sum_{k=1}^l \left(
  :\gamma^{+,k} (z) \psi^{-,k} (w): + w:\psi^{+,k}(-w) \gamma^{-,k}(-z):
  \right). $$

\begin{prop}
The above equations define a representation of $\widehat{\ospd}$
on $\F^{\otimes l}$ of central charge $l$.
\end{prop}

According to Feingold and Frenkel \cite{FF}, the Fourier
components of the generating functions
\begin{equation} \label{eq_fermi}
: \psi^{+,p}(z) \psi^{+,q} ( -z):, \, : \psi^{-,p}(z) \psi^{-,q} (
-z):, \,
 : \psi^{+,p}(z) \psi^{-,q} (z): + \hf \delta_{p,q}z^{-1}
\end{equation}
generate a representation of the twisted affine algebra $ {\rm
gl}^{(2)} ( 2l)$ of type $ A^{(2)}_{2l-1}$ on $\F^{\otimes l}$
with central charge $ 1$. On the other hand,  the Fourier
components of the following generating functions
\begin{equation} \label{eq_bose}
 -: \gamma^{+,p} (z) \gamma^{+,q} ( -z):, \quad
 - : \gamma^{-,p}(z) \gamma^{-,q} ( -z):, \quad
 - : \gamma^{+,p} (z) \gamma^{-,q} (z):
\end{equation}
generate a representation of the twisted affine algebra $ {\rm
gl}^{(2)} ( 2l)$ of type $ A^{(2)}_{2l-1}$ on $\F^{\otimes -l}$
with central charge $ -1$. We observe that the Fourier components
of the generating functions in \eqnref{eq_fermi} and
\eqnref{eq_bose} correspond to the same generators in $ {\rm
gl}^{(2)} ( 2l)$ with the appropriate coefficients as in
\eqnref{eq_fermi} and \eqnref{eq_bose} (cf. Eqs. (3.55--3.57),
pp139--140 in \cite{FF}; note that our convention here is a little
different). Then the diagonal action on $\F^{\otimes l}$ given by
\begin{eqnarray*}
  e^{pq}_{+} (z) & \equiv & \sum_{n \in \Z} e^{pq}_{+}(n)
                   z^{ -n -1}
  = : \psi^{+,p}(z) \psi^{+,q} ( -z):
    - : \gamma^{+,p} (z) \gamma^{+,q} ( -z):     \nonumber       \\
  e^{pq}_- (z) & \equiv & \sum_{n \in \Z} e^{pq}_-(n)
                   z^{ -n -1}
  = : \psi^{-,p}(z) \psi^{-,q} ( -z):
    - : \gamma^{-,p}(z) \gamma^{-,q} ( -z): \label{eq_series}       \\
  e^{pq} (z) & \equiv & \sum_{n \in \Z} e^{pq} (n)
                   z^{ -n -1}
  = : \psi^{+,p}(z) \psi^{-,q} (z):
    - : \gamma^{+,p} (z) \gamma^{-,q} (z):  + \hf  \delta_{p,q} z^{-1} \nonumber
\end{eqnarray*}
$ (p,q = 1, \dots, l  )$ is that of an affine algebra ${\rm
{gl}}^{(2)} (2l)$ of type $A^{(2)}_{2l -1}$ of central charge
zero. The horizontal subalgebra of the affine algebra ${\rm
{gl}}^{(2)} (2l)$ spanned by $e^{pq}_{+} (0), e_-^{pq} (0),$ $
e^{pq} (0) $ $(p, q = 1, \ldots, l)$ is isomorphic to the Lie
algebra ${ {\rm so}} (2l)$.

We claim that the action of the horizontal subalgebra ${ {\rm so}}
(2l)$ can be lifted to an action of the Lie group ${\rm Pin}(2l)$,
also cf. \cite{W1}. First let us take a digression to recall the
group ${\rm Pin}(2l)$.

The Pin group ${\rm Pin}(n)$ is the double covering group of $
{\rm O}(n)$, namely we have $$ 1 \longrightarrow \Z_2
\longrightarrow {\rm Pin}(n) \longrightarrow {\rm O}(n)
\longrightarrow 1. $$ We then define the Spin group ${\rm
Spin}(n)$ to be the inverse image of ${\rm SO} (n)$ under the
projection from ${\rm Pin} (n)$ to ${\rm O} (n) $. Then we have
the following exact sequence of Lie groups: $$ 1 \longrightarrow
\Z_2  \longrightarrow {\rm Spin} (n)
   \longrightarrow {\rm SO} (n) \longrightarrow 1.
$$
 Set $n = 2l$. Denote $ {\bf 1}_l = (1, 1, \ldots,1, 1)
\in \Z^l $ and ${\bar{\bf 1} }_l  = (1, 1, \ldots, 1, -1) \in
\Z^l$. The irreducible representations of ${\rm Spin} (2l)$ that
do not factor to ${\rm SO} (2l)$ are irreducible representations
of ${\rm {so}} (2l)$ of  highest weights parameterized by
\begin{eqnarray}
  \lambda = \hf {\bf 1}_l + (m_1, m_2, \ldots, m_l)    \label{eq_wtplus}
\end{eqnarray}
and
\begin{eqnarray}
\lambda = \hf {\bar{\bf 1} }_l + (m_1, m_2, \ldots, -m_l)
\label{eq_wtminus}
\end{eqnarray}
where $m_1 \geq \ldots \geq m_l \geq 0, m_i \in \Z.$ We are
interested in irreducible representations of ${\rm Pin} (2l)$
induced from irreducible representations of ${\rm Spin} (2l)$ with
highest weight of (\ref{eq_wtplus}) or (\ref{eq_wtminus}). When
restricted to ${\rm Spin} (2l)$, it will decompose into a sum of
the two irreducible representations of highest weights
(\ref{eq_wtplus}) and (\ref{eq_wtminus}). We will use $V({\rm
Pin}(2l), \la)$, where $\lambda = \hf |{\bf 1}_l | + (m_1, m_2,
\ldots, \overline{m}_l), m_l \geq 0$, to denote this irreducible
representation of ${\rm Pin} (2l)$. Denote by $$ \Sigma({\rm Pin})
= \{ \hf |{\bf 1}_l | + (m_1, m_2, \ldots, \overline{m}_l),
   m_1 \geq \ldots \geq m_l \geq 0, m_i \in \Z \}.  $$

The proof of the following lemma is straightforward.

\begin{lem}\label{ospdlem}
The action of ${\rm {so}}(2l)$ commutes with that of
$\widehat{\ospd}$ on $\F^{\otimes l}$.
\end{lem}

\begin{rem}
As a representation of ${\rm {so}}(2l)$, the Fock space
$\F^{\otimes l}$ is isomorphic to $\wedge (\C^l ) \wedge \wedge
(\C^{2l} \bigotimes \C^{\mathbb N}) \otimes S (\C^{2l} \bigotimes
\C^{\mathbb N})$, where $\wedge (\C^l )$ is the sum of two
half-spin representations and ${\rm {so}}(2l)$ acts on $\C^{2l}
\bigotimes \C^{\mathbb N}$ naturally by the left action on
$\C^{2l}$. The action of ${\rm {so}}(2l)$ can be lifted to ${\rm
Spin}(2l)$ which extends naturally to ${\rm Pin} (2l)$. It follows
that any irreducible representation of ${\rm Spin}(2l)$ appearing
in $\F^{\otimes l}$ cannot factor to ${\rm SO}(2l)$.  A similar
argument to the classical dual pair case \cite{H2} shows that
${\rm Pin} (2l)$ and $\hospd$ form a dual pair on $\F^{\otimes
l}$.
\end{rem}

The Cartan subalgebra $\hospd_0$ of $\hospd$ is spanned by the
basis elements $E_{aa}-E_{-a,-a}$, $a\in\hf\N$, and the central
element $1$. Let $\La\in(\hospd_0)^*$, the restricted dual of
$\hospd_0$.  We denote its components $\La(E_{aa}-E_{-a,-a})$ by
$\la_{a}$.

Given a Young diagram $\la=(\la_1,\la_2,\cdots,\la_l)$ of length
at most $l$ with
\begin{equation*}
\la_1\ge\la_2\ge\cdots\ge\la_{l-1}\ge\la_{l}\ge 0,
\end{equation*}
we may define a highest weight $\La(\la)\in(\hospd_0)^*$ as
follows:
\begin{eqnarray*}
\la_{j}&=& \langle\la_j'-j \rangle,\quad j\in\N,\\ \la_r&=&
\langle \la_{r+1/2}-(r-1/2) \rangle,\quad r\in\hf+\Z_+.
\end{eqnarray*}

For each $j\in\N$ define the matrix $ \tilde{X}^j$ to be the
matrix obtained from $ X^j$ by replacing its last column
$(\gamma_{-\hf}^{+,l},\cdots,\gamma_{-j+\hf}^{+,l},\psi_{-j}^{+,l},\cdots,\psi_{-j}^{+,l})$
by
\begin{equation*}
\left(-\gamma_{-\hf}^{+,l},\gamma_{-\frac{3}{2}}^{+,l},
-\gamma_{-\frac{5}{2}}^{+,l},\cdots,
(-1)^j\gamma_{-j+\hf}^{+,l},(-1)^{j+1}\psi_{-j}^{+,l},
\cdots,(-1)^{j+1}\psi_{-j}^{+,l} \right)
\end{equation*}
For $0\le r\le l$ let us, as usual, denote by $ \tilde{X}^j_r$ the
first $r\times r$ minor of $ \tilde{X}^j$.

\begin{thm}\label{hospdhwv}
\begin{enumerate}
\item As a ${\rm Pin}(2l)\times\hospd$-module, $\F^{\otimes l}$ is
completely reducible and decomposes into isotypic components as
follows:
\begin{equation*}
\F^{\otimes l}\cong\bigoplus_{\la \in\Sigma({\rm Pin})} V({\rm
Pin}(2l);\la) \bigotimes L(\hospd,\Lambda(\la)).
\end{equation*}

\item Let $\la=(\la_1,\la_2,\cdots,\la_l)$ be a Young diagram of
length at most $l$. The following vectors are joint highest weight
vectors of ${\rm so}(2l) \times \hospd$ in $\F^{\otimes l}$:
\begin{itemize}
\item[(a)]
${\rm det} X^1_{\la'_1}{\rm det} X^2_{\la'_2}\cdots{\rm det }
X^{\la_1}_{\la'_{\la_1}}$.
\item[(b)] $\psi_0^{-,l} {\rm det} \tilde{X}^1_{\la'_1}{\rm det} \tilde{X}^2_{\la'_2}
\cdots {\rm det} \tilde{X}^{\la_1}_{\la'_{\la_1}}$.
\end{itemize}
Furthermore the highest weight with respect to $\hospd$ of (a) and
(b) is $\La(\la)$, while their highest weights with respect to
${\rm so}(2l)$ are $\la=(\la_1+\hf,\la_2+\hf,\cdots,\la_l+\hf)$
and $\la=(\la_1+\hf,\la_2+\hf,\cdots,-\la_l-\hf)$, respectively.
\end{enumerate}
\end{thm}

\begin{proof}
From the second part we see that all irreducible modules of ${\rm
Pin}(2l)$  that do not factor through ${\rm O}(2l)$ appear in the
Fock space. The first part follows from this and the fact that
$\hospd$ and ${\rm Pin}(2l)$ form a dual pair.

In order to show that (a) is a highest weight vector it is enough
to check that it is annihilated by $e_+^{pq}(0)$, for all $p\not=
q$, due to \thmref{gltwohwv}, which is an easy computation.

The proof of (b) is somewhat more tedious, and we will omit the
details, as it is also a rather straightforward calculation. We
only remark that the crucial identity used in the computation is
again \lemref{identity}.

Finally, once the highest weight vectors are explicitly given, the
computation of weights of these vectors is straightforward.
\end{proof}


By composing the homomorphism $\hat{\phi}_{0}$ with the action of
$\hospd$ on $\F^{\otimes l}$, we realize a representation of
$\hsdzero$ on $\F^{\otimes l}$. Let us write down the action of
$\hsdzero$ on $\F^{\otimes l}$ explicitly. We introduce the
following generating functions:

\begin{eqnarray*}
{}^0 W_a^n (z) & =& \sum_{k \in\Z}  {}^0 W^n_{a, k} z^{-k -n -1},
a =0, 1, \\
{}^0 W^n_{\pm} (z) & =& \sum_{k \in \Z}  {}^0 W^n_{\pm, k} z^{-k
-n -1}.
\end{eqnarray*}

\begin{prop}\label{sdzeroff}
 The action of $\hsdzero$ on $\F^{\otimes l}$ is given in terms of
 the following generating functions:
\begin{eqnarray*}
   {}^0 W_0^n (z)
    & =& \sum_{p =1}^l
          \left(
                 :\psi^{+p} (z)\partial^n  \psi^{-p} (z):
                 -  z^{-1}:\partial^n (z\psi^{+p} (-z)) \psi^{-p} (-z):
          \right)                  \\
   {}^0 W_1^n (z)
    & =& -\sum_{p =1}^l
          \left(
                :\partial^n \gamma^{-p} (z) \gamma^{+p} (z):
                 + :\partial^n (\gamma^{+p} (-z)) \gamma^{-p} (-z):
          \right)                  \\
   {}^0 W_{\times}^n (z)
    & =& - \sum_{p =1}^l
          \left(
                :\partial^n \psi^{-p} (z) \gamma^{+p} (z):
                + :\partial^n (z \psi^{+p} (-z) ) \gamma^{-p} (-z):
          \right)
\end{eqnarray*}
\end{prop}

\begin{proof}
We will only give a proof of the first identity, as proofs for the
other two are obtained similarly.
  \begin{eqnarray*}
   {}^0 W_0^n (z)
        & =&\sum_{k \in \Z} t^k
           \left(
                 [D]_n + (- 1)^{k +1} [ -D -k]_n
           \right) M_0  z^{-k -n -1}   \\
        & =&  \sum_{k, j \in \Z}
           \left(
                 [-j]_n + (- 1)^{k +1} [ j -k]_n
           \right) E_{j -k, j} M_0  z^{-k -n -1}   \\
        & =&\sum_{k, j \in \Z}
           \left(
                 [-j]_n E_{j -k, j}
                 + (- 1)^{k +1} [-j]_n E_{-j, k -j}
           \right) M_0  z^{-k -n -1}   \\
       & =& \sum_{k, j \in \Z} \sum_{p =1}^l
           \left(
                 [-j]_n :\psi^{+p}_{k -j} \psi^{-p}_j :
                 + (- 1)^{k +1} [-j]_n: \psi^{+p}_j \psi^{-p}_{k -j}:
           \right) M_0  z^{-k -n -1}   \\
       & =& \sum_{p =1}^l
          \left(
                 :\psi^{+p} (z)\partial^n \psi^{-p} (z):
                 + \frac{1}{z} :\partial^n (z\psi^{+p} (-z)) \psi^{-p} (-z):
          \right).
\end{eqnarray*}
\end{proof}

\begin{rem} \label{remtaylor}
We recall that $E_0(z,w)+z^{-1}wE_0(-w, -z)$,
$E_1(z,w)+E_1(-w,-z)$ and $E_-(z,w)-wE_+(-w,-z)$ are the
generating functions for generators of $\hospd$. Indeed, one can
also express the field as
  $${}^0 W_0^n (z)=\partial_w^n (E_0(z,w)+z^{-1}wE_0(-w, -z))|_{w=z}.$$
 In other words, we have for $|w| >|z|$
 $$ E_0(z,w)+z^{-1}wE_0(-w, -z)
 =\sum_{n=0}^\infty \frac1{n!}{}^0 W_0^n (z) (w-z)^n.$$
 Similarly for $|w| >|z|$ we have
 \begin{eqnarray*}
  E_1(z,w)+E_1(-w,-z) =\sum_{n=0}^\infty \frac1{n!}{}^0 W_1^n (z)
  (w-z)^n, \\
  E_-(z,w)-wE_+(-w,-z) =\sum_{n=0}^\infty \frac1{n!}{}^0 W_{\times}^n (z) (w-z)^n.
  \end{eqnarray*}

\end{rem}

Now the irreducible module $L(\hospd, \Lambda(\la))$ pulls back to
via $\hat{\phi}_0$ a module over $\hsdzero$, which remains to be
irreducible thanks to \propref{propexactseq2} and
\propref{holoext}. For $\xi \in (\hsdzero)_0^*$, we denote by
$L(\hsdzero,\xi)$ the irreducible highest weight representation of
$\hsdzero$. Recall that the Cartan subalgebra of $\hsdzero$ is
spanned, aside from the central element, by
\begin{equation*}
\{D^n\parth\theta|n\ {\rm odd}\}\cup\{(D+\hf)^n\theta\parth|n\
{\rm odd}\},
\end{equation*}
so that its generating functions are ${}^{0}\Delta_0 (x)=-\xi
(\parth \theta \sinh (xD) )$ and ${}^0\Delta_1 (x)= -\xi (
\theta\parth \sinh (x( D + 1/2) )$.

\begin{prop}\label{hsdzeroprop1}
Let $L(\hospd,\Lambda)$ be the irreducible highest weight
representation of highest weight $\Lambda$ such that
$\Lambda(E_{aa}-E_{-a,-a}) = \lambda_a$, $a\in\hf\N$. Then via
$\hat{\phi}_0$ the module $L(\hgltwo,\La)$ pulls back to the
irreducible highest weight module $L(\hsdzero,\La_{\hsdzero})$,
where the highest weight $\La_{\hsdzero}$ is specified by the
generating functions
\begin{align*}
{}^0\Delta_0(x)&= -\sum_{i\in\hf+ \Z_+} \la_{rr}{\rm sinh}(-rx),\\
{}^0\Delta_1(x)&= -\sum_{i\in\N} \la_{ii}{\rm sinh}(-ix).
\end{align*}
\end{prop}

\begin{proof}
The generating functions for $\hsdzero$ can be computed using
\propref{propext2}.  Irreducibility is a consequence of
\propref{propexactseq2} and \propref{holoext}.
\end{proof}

Now via $\hat{\phi}_s$, $s\not=\hf\Z$, the irreducible
representation $L(\hgltwo, \Lambda(\la))$ pulls back to an
irreducible representation of $\hsdzero$ by
\propref{propexactseq2}.  Again using \propref{propext2} we can
show the following.

\begin{prop}\label{hsdzeroprop2}
Let $L(\hgltwo,\Lambda)$ be the irreducible highest weight
representation of highest weight $\Lambda$ and central charge $c$
such that $\Lambda(E_{aa}) = \lambda_a$. Then via $\hat{\phi}_s$,
$s\not=\hf\Z$, the module $L(\hgltwo,\La)$ pulls back to the
irreducible highest weight module $L(\hsdzero,\La_{\hsdzero})$,
where the highest weight $\La_{\hsdzero}$ is specified by the
generating functions
\begin{align*}
{}^0\Delta_0(x)&=-\sum_{r\in\hf+\Z}{\rm
sinh}(-rx)\la_{rr}-\frac{{\rm cosh}(sx)-1}{{\rm sinh}(x/2)}c,\\
{}^0\Delta_1(x)&=-\sum_{i\in\Z}{\rm sinh}((-i+s)x)\la_{ii}
+\frac{{\rm cosh}((s-\hf)x)-{\rm cosh}(x/2)}{{\rm sinh}(x/2)}c.
\end{align*}
\end{prop}

\begin{rem}
It follows that $\hsdzero$ and ${\rm Pin}(2l)$ form a dual pair in
the sense of Howe on $\F^{\otimes l}$.  The decomposition of
$\F^{\otimes l}$ with respect to the joint action
$\hsdzero\times{\rm Pin}(2l)$ is easily obtained from
\thmref{hospdhwv} and \propref{hsdzeroprop1}. We further note that
\thmref{gltwohwv} together with \propref{hsdzeroprop2} imply that
${\hsdzero}$ and ${\rm GL}(l)$ form a dual pair on $\F^{\otimes
l}$ as well.
\end{rem}

\subsection{The case of $\hpd_{\pm\pm}$ and ${}^{\pm\pm}\widehat{\sd}$}

We may also compose the homomorphism $\hat{\phi}_{0}$ with the
action of $\hpd_{\pm\pm}$ on $\F^{\otimes l}$.  This way we obtain
a representation of ${}^{\pm\pm}\hsd$ on $\F^{\otimes l}$ which we
will describe explicitly. We introduce the following generating
functions (we recall that the elements ${}^{\pm\pm}
W^n_k,{}^{\pm\pm} W^n_{+,k}$ and ${}^{\pm\pm} W^n_{-,k}$ in
${}^{\pm\pm}\hsd$ were introduced at the end of
section~\ref{sect_subalg}):

\begin{eqnarray*}
  {}^{\pm\pm} W^n (z)
   & =& \sum_{k \in\Z}  {}^{\pm\pm} W^n_k z^{-k -n -1}, \\
  {}^{\pm\pm} W^n_{\pm} (z)
   & =& \sum_{k\in\Z}  {}^{\pm\pm} W^n_{+, k} z^{-k -n -1}, \\
  {}^{\pm\pm} W^n_{\pm} (z)
   & =& \sum_{k\in\Z}  {}^{\pm\pm} W^n_{-, k} z^{-k -n -1}.
\end{eqnarray*}

\begin{prop}
 The action of ${}^{+\pm}\hsd$ on $\F^{\otimes l}$ is given in terms of
 the following generating functions:
 \begin{eqnarray*}
   {}^{+\pm} W^n (z)
    & =& \sum_{p =1}^l
          \left(
                 :\psi^{+p} (z)\partial^n \psi^{-p} (z):
                 + :\partial^n \gamma^{+p} (z) \gamma^{-p} (z):
          \right)                  \\
   {}^{+\pm} W_+^n (z)
    & =&-\sum_{p =1}^l
          \left(
                :\gamma^{+p} (z)\partial^n \psi^{-p} (z):
                 \mp :\partial^n (\gamma^{+p} (z)) \psi^{-p} (z):
          \right)                  \\
   {}^{+\pm} W_-^n (z)
    & =&  \sum_{p =1}^l
          \left(
                :\psi^{+p} (z)\partial^n \gamma^{-p} (z):
                 \pm :\partial^n \psi^{+p} (z) \psi^{-p} (z):
          \right)
 \end{eqnarray*}
\end{prop}

\begin{proof}
We will only give a proof of the second identity, as proofs for
the other two are obtained similarly.
  \begin{eqnarray*}
   {}^{+\pm} W_+^n (z)
        & =&\sum_{k \in \Z} t^k
           \left(
                 [D]_n \mp [ -D -k-1]_n
           \right) M_-  z^{-k -n -1}   \\
        & =&  \sum_{k, j \in \Z}
           \left(
                 [-j]_n \mp [j -k-1]_n
           \right) E_{j -k, j} M_-  z^{-k -n -1}   \\
        & =&  \sum_{k, j \in \Z}
            -\left(
                 [-j]_n \mp [j -k-1]_n
           \right) \sum_{p =1}^l\gamma^{+p}_{k-j+1/2} \psi^{-p}_j z^{-k -n -1}   \\
       & =& -\sum_{p =1}^l
          \left(
                 :\gamma^{+p} (z)\partial^n \psi^{-p} (z):
                 \mp :\partial^n (\gamma^{+p} (z)) \psi^{-p} (z):
          \right).
\end{eqnarray*}
\end{proof}

\begin{rem}
Let ${\mathfrak P}^{\pm} = \{ {}^{+\pm} W^n_k,{}^{+\pm}
W^n_{+,k},{}^{+\pm} W^n_{-,k}  \mid n +k \geq 0, k \in \Z, n \in
\Z_+\}$. One can show that $\widehat{\mathfrak P}^{\pm} =
{\mathfrak P}^{\pm} \bigoplus \C C $ is a parabolic subalgebra of
${}^{+\pm}\hsd$. As in the $\hsd$ case, we define the vacuum
${}^{+\pm}\hsd$-module $M_c({}^{+\pm}\hsd)$ ($c\in\C$) $$
  M_c({}^{+\pm}\hsd) = {\mathcal U} ({}^{+\pm}\hsd)
                 \bigotimes_{ {\mathcal U}(\widehat{\mathfrak P}^{\pm})}
                 \C_c
$$ and its irreducible quotient $V_c({}^{+\pm}\hsd)$. Then
$M_c({}^{+\pm}\hsd)$ and $V_c({}^{+\pm}\hsd)$ also carries vertex
superalgebra structures.
\end{rem}

\begin{prop}
 The action of \;${}^{-\pm}\hsd$ on $\F^{\otimes l}$ is given in terms of
 the following generating functions:
 \begin{eqnarray*}
   {}^{-\pm} W^n (z)
    & =& \sum_{p =1}^l
          \left(
                 :\psi^{+p} (z)\partial^n \psi^{-p} (z):
                 - :\partial^n \gamma^{+p} (-z) \gamma^{-p} (-z):
          \right)                  \\
   {}^{-\pm} W_+^n (z)
    & =&\sum_{p =1}^l
          \left(
                :\gamma^{+p} (z)\partial^n \psi^{-p} (z):
                 \pm :\partial^n (\gamma^{+p} (-z)) \psi^{-p} (-z):
          \right)                  \\
   {}^{-\pm} W_-^n (z)
    & =&  \sum_{p =1}^l
          \left(
                :\psi^{+p} (z)\partial^n \gamma^{-p} (z):
                 \mp :\partial^n (\psi^{+p} (z)) \gamma^{-p} (z):
          \right)
 \end{eqnarray*}
\end{prop}

\begin{proof}
We will give a proof of the second identity only.
  \begin{eqnarray*}
   {}^{-\pm} W_+^n (z)
        & =&\sum_{k \in \Z} t^k
           \left(
                 [D]_n \pm (-1)^{k+1} [ -D -k-1]_n
           \right) M_-  z^{-k -n -1}   \\
        & =&  \sum_{k, j \in \Z}
           \left(
                 [-j]_n \pm (-1)^{k+1} [j -k-1]_n
           \right) E_{j-k, j} M_-  z^{-k -n -1}   \\
        & =&  \sum_{k, j \in \Z}
            \left(
                 [-j]_n \pm (-1)^{k+1} [j -k-1]_n
           \right) \sum_{p =1}^l\gamma^{+p}_{k-j+1/2} \psi^{-p}_{j} z^{-k -n -1}   \\
       & =& \sum_{p =1}^l
          \left(
                 :\gamma^{+p} (z)\partial^n \psi^{-p} (z):
                 \pm :\partial^n (\gamma^{+p} (-z)) \psi^{-p} (-z):
          \right).
\end{eqnarray*}

\end{proof}

\begin{rem} The Fock space $\F^{\otimes l}$ is not completely reducible
with respect to the action of $\hpd_{\pm\pm}$. Similar remarks as
\remref{remtaylor} hold in the cases of ${}^{\pm\pm}\hsd$.
\end{rem}
\section{Conclusion and discussion}

In this paper we have classified all anti-involutions of Lie
superalgebra $\hsd$ preserving the principal gradation, where
$\hsd$ is the central extension of the Lie superalgebra of
differential operators on the super circle $S^{1|1}$. There are
five familes of them, while the anti-involutions within a family
is related to each other by a spectral flow. We found close
relations between the subalgebras ${}^0\hsd, {}^{\pm\pm}\hsd$
fixed by the five anti-involutions and subalgebras $\hospd$,
$\hpd_{\pm\pm}$ of $\widehat{gl}_{\infty|\infty}$. We realize
these Lie superalgebras by a free field realization, and further
establish dualities between them and certain finite-dimensional
classical Lie groups on Fock spaces. Our construction generalizes
\cite{KWY}.

There are various interesting questions arising from our current
work. Below we will list some of them. Some of these questions
have been better understood in the $\WW_{1+\infty}$ case.

\begin{enumerate}
\item The Lie superalgebra $\hsd$, i.e. the super
extension of $\W_{1+\infty}$ was first studied by Manin and Radul
\cite{MR} in the context of sypersymmetric KP hierarchy. It is
natural to ask to construct the hierarchy corresponding to the Lie
superalgebras $\hsdzero$ and ${}^{\pm\pm}\hsd$.


\item Calculate the principle $q$-character formulas
for the irreducible modules of $\hsd$, $\hsdzero$ and
${}^{\pm\pm}\hsd$ that appear in the free field realizations.

\item Identify the vertex superalgebras $V_l(\hsd)$ and
$V_l({}^{+\pm}\hsd)$ of central charge $l \in \N$ with more
familiar $\mathcal W$ superalgebras (cf. \cite{BS}). Such a
question has been answered for ${\mathcal W}_{1+\infty}$ when $l
\in \N$ \cite{FKRW} and when $l =-1$ \cite{W2}.

\item Classify the irreducible modules of the vertex
superalgebras $V_l(\hsd)$ and $V_l({}^{+\pm}\hsd)$.

\item What is the fusion ring of the vertex
superalgebras $V_l(\hsd)$ and $V_l({}^{+\pm}\hsd)$ of central
charge $l\in \N$? The duality between $(GL(l), {\mathcal
W}_{1+\infty}|_{c=l})$ lead the authors in \cite{FKRW} to
conjecture that the fusion ring for ${\mathcal
W}_{1+\infty}|_{c=l}$ is isomorphic to the representation ring of
$GL(l)$. (Indeed one can argue that the former contains the
later). It is natural to conjecture that the fusion ring for
$\hsd$ of central charge $l$ is also isomorphic to the
representation ring of $GL(l)$ based on our duality results.

\item We have studied the superalgebra of differential
operator on the super circle $S^{1|N}$ with $N=1$ extended
symmetry. One may ask similar questions for the super circle
$S^{1|N}$ for a general $N$.
\end{enumerate}

\end{document}